\documentclass[12pt,reqno]{amsart}
\usepackage{amsmath}
\usepackage{amssymb}
\usepackage{verbatim}
\usepackage{mathrsfs}

\usepackage[top=0.6in,bottom=0.6in,left=0.7in,right=0.7in]{geometry}

\newcommand{\C}{\mathbb{C}}
\newcommand{\N}{\mathbb{N}}    
\newcommand{\R}{\mathbb{R}}    
\newcommand{\Z}{\mathbb{Z}}

\newcommand{\dR}{\mathbb{R}^d}

\newcommand{\dZ}{\mathbb{Z}^d}

\newcommand{\be}{ \begin{equation} }
\newcommand{\ee}{ \end{equation} }
\newcommand{\bp}{ \begin{proof} }
\newcommand{\ep}{\hfill  \end{proof} }
\newcommand{\gl}{\lambda}
\newcommand{\gep}{\varepsilon}
\newcommand{\dm}{\mathsf{M}}
\newcommand{\dmfc}{\Omega_{\dm}} 
\newcommand{\si}{\mathtt{S}}            
\newcommand{\supp}{\operatorname{supp}}
\newcommand{\diag}{\operatorname{diag}}
\newcommand{\td}{\boldsymbol{\delta}}
\newcommand{\vertiii}[1]{{\left\vert\kern-0.45ex\left\vert\kern-0.45ex\left\vert #1\right\vert\kern-0.45ex\right\vert\kern-0.45ex\right\vert}}  
\newcommand{\mnorm}[1]{\vertiii{#1}}
\newcommand{\imply}{ \Longrightarrow }
\newcommand{\AS}{\operatorname{\mathsf{AS}}}
\newcommand{\setsp}{\;:\;}     
\newcommand{\tp}{\mathsf{T}}  
\newcommand{\len}{\operatorname{len}}
\newcommand{\Lp}[1]{L_{#1}(\mathbb{R})}
\newcommand{\dLp}[1]{L_{#1}(\mathbb{R}^d)}
\newcommand{\dTLp}[1]{L_{#1}(\mathbb{T}^d)}
\newcommand{\V}{\mathscr{V}}   
\newcommand{\wh}{\widehat}
\renewcommand{\le}{\leqslant}
\renewcommand{\ge}{\geqslant}
\newcommand{\bs}{\backslash}
\newcommand{\ol}{\overline}
\newcommand{\la}{\langle}
\newcommand{\ra}{\rangle}
\newcommand{\fr}{\mathcal{F}}
\newcommand{\bo}{\mathscr{O}} 

\newtheorem{lemma}{Lemma}
\newtheorem{prop}[lemma]{Proposition}
\newtheorem{cor}[lemma]{Corollary}
\newtheorem{theorem}[lemma]{Theorem}

\numberwithin{equation}{section}
\numberwithin{lemma}{section}

\begin{document}

\title{Homogeneous Wavelets and Framelets with the Refinable Structure}

\author{Bin Han}

\thanks{Research supported in part by NSERC Canada under Grant
RGP 228051. 
}

\address{Department of Mathematical and Statistical Sciences,
University of Alberta, Edmonton,\quad Alberta, Canada T6G 2G1.
\quad {\tt bhan@ualberta.ca}\quad
{\tt http://www.ualberta.ca/$\sim$bhan}
}

\dedicatory{\centerline{Dedicated to the memory of Professor Minde Cheng}}

\makeatletter \@addtoreset{equation}{section} \makeatother
\begin{abstract}
Homogeneous wavelets and framelets have been extensively investigated in the classical theory of wavelets and they are often constructed from refinable functions via the multiresolution analysis. On the other hand, nonhomogeneous wavelets and framelets enjoy many desirable theoretical properties and are often intrinsically linked to the refinable structure and multiresolution analysis. In this paper we shall provide a comprehensive study on connecting
homogeneous wavelets and framelets to nonhomogeneous ones with the refinable structure. This allows us to understand better the structure of homogeneous wavelets and framelets as well as their connections to the refinable structure and multiresolution analysis.
\end{abstract}

\keywords{Homogeneous wavelets and framelets, nonhomogeneous wavelets and framelets, refinable structure, shift-invariant spaces, multiresolution analysis, Schur decomposition for Hermite matrices of measurable functions, singular value decomposition for matrices of measurable functions}

\subjclass[2010]{42A40, 42C15, 41A30} \maketitle

\pagenumbering{arabic}

\section{Introduction and Motivations}

Since the birth of wavelet theory in 1980's, wavelets and framelets have been extensively investigated and successfully used in many applications (\cite{chubook92,daubook92,malbook09,meybook90}). In particular, homogeneous wavelets and framelets are the major topics in wavelet theory and they are often constructed from refinable (vector) functions via the multiresolution analysis and the refinable structure. To better understand homogeneous wavelets and framelets,
in this paper we shall link homogeneous wavelets and framelets to nonhomogeneous ones with the refinable structure.

For a function $f: \dR\rightarrow \C$ and a $d\times d$ real-valued matrix $U$, throughout this paper we shall adopt the following notation:
\be \label{f:ds}
f_{U;k}(x):=|\det(U)|^{1/2} f(Ux-k),\qquad x,k\in \dR.
\ee
For an integrable function $f\in \dLp{1}$, its Fourier transform $\wh{f}$ in this paper is defined to be $\wh{f}(\xi):=\int_{\dR} f(x) e^{-ix\cdot \xi} dx$ for $\xi\in \dR$.
The Fourier transform can be naturally extended to square integrable functions in $\dLp{2}$.

For a $d\times d$ real-valued invertible matrix $\dm$ and a countable subset $\Psi$ of square integrable functions in $\dLp{2}$, we recall that \emph{a homogeneous $\dm$-affine system} $\AS(\Psi)$ is defined to be
\be \label{has}
\AS(\Psi):=\{\psi_{\dm^j;k} \setsp j\in \Z, k\in \dZ, \psi\in \Psi\}.
\ee
A homogeneous affine system $\AS(\Psi)$ is often regarded as the natural discretization of a continuous wavelet transform using the generating wavelet $\Psi$ (see \cite{chubook92,daubook92}) and has been widely studied in wavelet analysis and applied harmonic analysis (\cite{chubook92,daubook92,longbook95,malbook09,meybook90} and many references therein).
The generating set $\Psi$ is called \emph{a homogeneous orthogonal $\dm$-wavelet} in $\dLp{2}$ if its homogeneous $\dm$-affine system $\AS(\Psi)$ is an orthonormal basis of $\dLp{2}$.
One-dimensional dyadic (i.e., $\dm=2$) homogeneous orthogonal $2$-wavelets $\Psi$, with $\Psi$ being a singleton,  have been extensively and systematically studied in the books of Chui \cite{chubook92}, Daubechies~\cite{daubook92},
Hern\'andez and Weiss~\cite{hwbook96},
Meyer~\cite{meybook90}, Mallat~\cite{malbook09}, and many other books and numerous papers in the literature.

We say that a $d\times d$ real-valued matrix $\dm$ is \emph{expansive} if all its eigenvalues are greater than one in modulus. The commonly used expansive integer matrices in the literature of wavelet theory include $\dm=2I_d$ (the dyadic diagonal matrix) and the two-dimensional matrices
\[
\dm_{\sqrt{2}}=\left[ \begin{matrix} 1 &1\\ 1 &-1\end{matrix}\right],\qquad
\dm_{\sqrt{3}}=\left[ \begin{matrix} 1 &-2\\ 2 &-1\end{matrix}\right].
\]
As discussed in \cite{han12,han10},
a closely related notion to a homogeneous affine system is a nonhomogeneous affine system. Let $\Phi$ and $\Psi$ be countable subsets of $\dLp{2}$. For every integer $J\in \Z$, \emph{a nonhomogeneous $\dm$-affine system} $\AS_J(\Phi;\Psi)$ is defined to be
\be \label{as}
\AS_J(\Phi;\Psi):=\{\phi_{\dm^J;k} \setsp k\in \dZ,\phi\in \Phi\}\cup
\{\psi_{\dm^j;k} \setsp j\ge J, k\in \dZ, \psi\in \Psi\}.
\ee
By a simple scaling argument, it is trivial to observe that if $\AS_J(\Phi;\Psi)$ is an orthonormal basis of $\dLp{2}$ for some integer $J\in \Z$, then it is an orthonormal basis of $\dLp{2}$ for every $J\in \Z$.
A generating set $\{\Phi;\Psi\}$ in $\dLp{2}$ is called \emph{an orthogonal $\dm$-wavelet} in $\dLp{2}$ if its nonhomogeneous $\dm$-affine system $\AS_0(\Phi;\Psi)$ is an orthonormal basis of $\dLp{2}$.
Under the condition that $\dm$ is a $d\times d$ expansive real-valued matrix and $\sum_{\phi\in \Phi} \|\phi\|_{\dLp{2}}^2<\infty$, if $\{\Phi;\Psi\}$ is an orthogonal $\dm$-wavelet in $\dLp{2}$, then we know from \cite[Theorem~6]{han12} that $\Psi$ must be a homogeneous orthogonal $\dm$-wavelet in $\dLp{2}$. Moreover, suppose that both $\Phi=\{\phi^1,\ldots,\phi^r\}$ and $\Psi=\{\psi^1,\ldots,\psi^s\}$ are finite subsets of $\dLp{2}$ with $r,s\in \N$.
Define two vector functions $\phi\in (\dLp{2})^r$ and $\psi\in (\dLp{2})^s$ by
$\phi:=(\phi^1,\ldots,\phi^r)^\tp$ and
$\psi:=(\psi^1,\ldots,\psi^s)^\tp$.
It has been shown in \cite[Theorem~7]{han12} that if
$\dm$ is a $d\times d$ invertible integer matrix and if $\{\Phi;\Psi\}$ is an orthogonal $\dm$-wavelet in $\dLp{2}$, then $s=r(|\det(\dm)|-1)$ and $\{\Phi;\Psi\}$ must possess the following refinable structure:
\be \label{refstr}
\wh{\phi}(\dm^\tp\xi)=\wh{a}(\xi)\wh{\phi}(\xi)
\qquad\mbox{and}\qquad
\wh{\psi}(\dm^\tp\xi)=\wh{b}(\xi)\wh{\phi}(\xi),\qquad a.e.\; \xi\in \dR,
\ee
where $\wh{a}$ and $\wh{b}$ are $r\times r$ and $s\times r$ matrices of $2\pi\dZ$-periodic measurable functions in $\dTLp{2}$, respectively, such that the filter bank $\{\wh{a};\wh{b}\}$ must form \emph{an orthogonal $\dm$-wavelet filter bank}, i.e.,
\begin{align}
&\ol{\wh{a}(\xi)}^\tp\wh{a}(\xi)+
\ol{\wh{b}(\xi)}^\tp\wh{b}(\xi)=I_r,\qquad a.e.\, \xi\in \dR, \label{owfb1}\\
&\ol{\wh{a}(\xi)}^\tp\wh{a}(\xi+2\pi \omega)+
\ol{\wh{b}(\xi)}^\tp\wh{b}(\xi+2\pi \omega)=0,\qquad a.e.\, \xi\in \dR, \label{owfb0}
\end{align}
for all $\omega\in \dmfc\bs\{0\}$,
where $\dmfc:=[(\dm^\tp)^{-1}\dZ]\cap [0,1)^d$.
That is, nonhomogeneous orthogonal wavelets are intrinsically connected to the refinable structure and filter banks  in \eqref{refstr}. Therefore, it is of interest to link a homogeneous wavelet with a nonhomogeneous wavelet.

A frame is a generalization of a basis in a Hilbert space. The notion of a frame was first introduced in Duffin and Schaeffer \cite{ds52} in the setting of nonharmonic Fourier series.
For a countable subset $\Psi\subset \dLp{2}$, we say that $\Psi$ is \emph{a homogeneous $\dm$-framelet} in $\dLp{2}$ if its homogeneous $\dm$-affine system $\AS(\Psi)$ is a frame for $\dLp{2}$, that is, there exist positive constants $C_1$ and $C_2$ such that
\be \label{hframelet}
C_1\|f\|^2_{\dLp{2}}\le
\sum_{j\in \Z} \sum_{\psi\in \Psi}\sum_{k\in \dZ} |\la f, \psi_{\dm^j;k}\ra|^2\le C_2\|f\|_{\dLp{2}}^2,\qquad \forall\, f\in \dLp{2},
\ee
where the optimal constants $C_1$ and $C_2$ in \eqref{hframelet} are called the lower and upper frame bounds, respectively.
In particular, the generating set $\Psi$ is called \emph{a homogeneous tight $\dm$-framelet} if \eqref{hframelet} holds with $C_1=C_2=1$, i.e., $\AS(\Psi)$ is a (normalized) tight frame for $\dLp{2}$. A homogeneous orthogonal $\dm$-wavelet in $\dLp{2}$ is obviously a homogeneous tight $\dm$-framelet in $\dLp{2}$. By a straightforward argument, it is well known that $\Psi$ is a homogeneous orthogonal $\dm$-wavelet in $\dLp{2}$ if and only if $\Psi$ is a homogeneous tight $\dm$-framelet in $\dLp{2}$ such that $\|\psi\|_{\dLp{2}}=1$ for all $\psi\in \Psi$.
The recent revived interest to study frames is largely inspired by the pioneering work of Daubechies, Grossmann and Meyer \cite{dgm86}, where one-dimensional bandlimited homogeneous tight $2$-framelets in $\Lp{2}$ have been constructed. Theory and construction of homogeneous framelets in $\dLp{2}$ have been extensively studied in numerous papers, e.g. see
\cite{bow00,chs02,css98,dh04,dhrs03,hanbook,han15,han14,han13,han12,han10,han03,han97,hanmsc94,hm05,hm03,hz14,hwbook96,longbook95,rs97df,rs97} and many references therein. Let $\Phi$ and $\Psi$ be subsets of $\dLp{2}$. By a simple scaling argument, it is known in \cite[Proposition~4]{han12} that if $\AS_J(\Phi;\Psi)$ is a frame of $\dLp{2}$ for some integer $J\in \Z$, then $\AS_J(\Phi;\Psi)$ is a frame of $\dLp{2}$ for every $J\in \Z$ with the same lower and upper frame bounds. Therefore, we call $\{\Phi;\Psi\}$ \emph{an $\dm$-framelet} in $\dLp{2}$ if $\AS_0(\Phi;\Psi)$ is a frame for $\dLp{2}$. If in addition $\dm$ is expansive and $\sum_{\phi\in \Phi} \|\phi\|_{\dLp{2}}^2<\infty$, then \cite[Proposition~4]{han12} tells us that $\Psi$ must be a homogeneous framelet in $\dLp{2}$. In particular, $\{\Phi;\Psi\}$ is called \emph{a tight $\dm$-framelet} in $\dLp{2}$ if $\AS_0(\Phi;\Psi)$ is a (normalized) tight frame for $\dLp{2}$.
It has been shown in \cite{hanrefstr} (also see Theorem~\ref{thm:ntf} in this paper) that a nonhomogeneous tight framelet is intrinsically linked to the refinable structure in \eqref{refstr}.

For a homogeneous $\dm$-framelet $\Psi$ in $\dLp{2}$, the $\dm$-affine system $\AS(\Psi)$ is a frame in $\dLp{2}$ and its associated frame operator $\fr: \dLp{2}\rightarrow \dLp{2}$ given by
\[
\fr (f):=\sum_{h\in \AS(\Psi)} \la f, h\ra h,\qquad f\in \dLp{2}
\]
must be a bijective bounded linear operator on $\dLp{2}$ with the bounded inverse operator $\fr^{-1}$. Consequently, it is well known that $\{ \fr^{-1} h\}_{h\in \AS(\Psi)}$ must be a dual frame of the given frame $\AS(\Psi)$ (called the canonical dual frame of $\AS(\Psi)$, e.g., see \cite{daubook92}), i.e., $\{ \fr^{-1} h\}_{h\in \AS(\Psi)}$ must be a frame for $\dLp{2}$ and
$\la f,g\ra=\sum_{h\in \AS(\Psi)} \la f, \fr^{-1} h\ra \la h, g\ra$ for all $f,g\in \dLp{2}$
with the series  converging absolutely.
However, the dual frame $\{ \fr^{-1} h\}_{h\in \AS(\Psi)}$ may not have the affine structure in \eqref{has} (e.g., see \cite{daubook92,dh02}), i.e., there may not exist $\tilde{\Psi}\subset \dLp{2}$ such that $\{ \fr^{-1} h\}_{h\in \AS(\Psi)}=\AS(\tilde{\Psi})$. Besides its canonical dual frame, a frame may have many other dual frames with or without the affine structure in \eqref{has}.
This naturally leads to the notion of homogeneous dual framelets.
Let $\Psi$ and $\tilde{\Psi}$ be subsets of $\dLp{2}$ such that $\Psi$ and $\tilde{\Psi}$ have the same cardinality. Throughout this paper we shall use the notation $\sim: \Psi\rightarrow \tilde{\Psi}$ with $\psi\mapsto \tilde{\psi}$ to stand for a given bijection between $\Psi$ and $\tilde{\Psi}$, that is, $(\psi,\tilde{\psi})$ is always regarded as a pair with $\psi\in \Psi$ and $\tilde{\psi}\in \tilde{\Psi}$.
We say that $(\tilde{\Psi},\Psi)$ is \emph{a homogeneous dual $\dm$-framelet} in $\dLp{2}$ if (i) each of $\Psi$ and $\tilde{\Psi}$ is an $\dm$-framelet in $\dLp{2}$, i.e., both $\AS(\Psi)$ and $\AS(\tilde{\Psi})$ are frames for $\dLp{2}$, and (ii) the following identity holds
\be \label{hdf}
\sum_{j\in \Z} \sum_{\psi\in \Psi} \sum_{k\in \dZ} \la f, \psi_{\dm^j;k}\ra\la \tilde{\psi}_{\dm^j;k},g\ra=\la f,g\ra,\qquad \forall\, f,g\in \dLp{2}
\ee
with the series converging absolutely.

A homogeneous $\dm$-framelet is also related to the notion of a homogeneous $\dm$-wavelet in $\dLp{2}$.
Let us recall the definition of a homogeneous Riesz wavelet in $\dLp{2}$. Let $\Psi$ be a countable subset of $\dLp{2}$. We say that $\Psi$ is \emph{a homogeneous Riesz $\dm$-wavelet} in $\dLp{2}$ if $\AS(\Psi)$ is a Riesz basis of $\dLp{2}$, that is, the linear span of $\AS(\Psi)$ is dense in $\dLp{2}$ and $\AS(\Psi)$ is a Riesz sequence in $\dLp{2}$, i.e.,
there exist positive constants $C_3$ and $C_4$ such that
\be \label{rieszwavelet}
C_3\sum_{j\in \Z} \sum_{\psi\in \Psi} \sum_{k\in \dZ}
|w_{j;k,\psi}|^2\le \Big\|
\sum_{j\in \Z} \sum_{\psi\in \Psi}
\sum_{k\in \dZ} w_{j;k,\psi} \psi_{\dm^j;k}\Big\|_{\dLp{2}}^2\le
C_4\sum_{j\in \Z} \sum_{\psi\in \Psi} \sum_{k\in \dZ}
|w_{j;k,\psi}|^2
\ee
for all finitely supported sequences $\{ w_{j;k,\psi} \}_{j\in \Z, k\in \dZ, \psi\in \Psi}$. The constant $C_3$ is called a lower Riesz bound and $C_4$ is an upper Riesz bound of $\AS(\Psi)$.
It is well known that a Riesz basis in $\dLp{2}$ is simply a frame in $\dLp{2}$ which is $l_2$-linearly independent, i.e., if $\sum_{h\in \AS(\Psi)} c_h h=0$ with $\sum_{h\in \AS(\Psi)} |c_h|^2<\infty$, then $c_h=0$ for all $h\in \AS(\Psi)$.
In other words, a Riesz basis is just a frame without redundancy.

For subsets $\Phi,\Psi\subset \dLp{2}$, we say that $\{\Phi;\Psi\}$ is \emph{a Riesz $\dm$-wavelet} in $\dLp{2}$ if $\AS_0(\Phi;\Psi)$ is a Riesz basis of $\dLp{2}$. It is known in \cite[Theorem~6]{han12} that if $\AS_J(\Phi;\Psi)$ is a Riesz basis of $\dLp{2}$ for some integer $J$, then $\AS_J(\Phi;\Psi)$ is a Riesz basis for all $J\in \Z$. Moreover, if $\sum_{\phi\in \Phi}\|\phi\|_{\dLp{2}}^2<\infty$ and $\dm$ is expansive, then a Riesz $\dm$-wavelet $\{\Phi;\Psi\}$ in $\dLp{2}$ always leads to a homogeneous Riesz $\dm$-wavelet $\Psi$ in $\dLp{2}$.

Nonhomogeneous affine systems have been comprehensively studied in Han~\cite{han12,han10}. As demonstrated in \cite{han12,han10}, in contrast to the study of homogeneous wavelets and framelets (e.g., see \cite{bow00,han97,hanmsc94,hwbook96,longbook95,rs97}), it is often much simpler to characterize and construct nonhomogeneous wavelets and framelets than homogeneous ones.
Due to many desirable properties of nonhomogeneous affine systems, it is of interest in both theory and application to investigate when a homogeneous affine system can be linked to a nonhomogeneous affine system.
Along this direction, \cite{ams14,aps16} and \cite[Section~4.5]{hanbook} have studied the connections of a homogeneous tight/dual framelet to a nonhomogeneous tight/dual framelet, while the connection of a homogeneous wavelet with a multiresolution analysis has been discussed in \cite{bow03,hanbook,han12,han95,hanmsc94,kkl01,ltw01,zal99} and references therein.
In this paper we shall comprehensively study the connections of homogeneous wavelets and framelets with nonhomogeneous wavelets and framelets and their connections to the refinable structure.
We shall obtain almost complete satisfactory answers to this topic.
In particular, our results on homogeneous framelets include all the results in \cite{ams14,aps16,hanbook} as special cases.

The structure of the paper is as follows.
In Section~2 we shall recall some known results on shift-invariant subspaces of $\dLp{2}$. Then we shall discuss the Schur decomposition for Hermite matrices of measurable functions and singular value decomposition for general matrices of measurable functions.
In Section~3, we shall study the connections of homogeneous framelets and homogeneous tight framelets with nonhomogeneous ones and the refinable structure. In particular, the intrinsic relation between nonhomogeneous tight framelets and the refinable structure will be presented in Section~3.
In Section~4, we link homogeneous dual framelets to nonhomogeneous dual framelets. Comparing with \cite{ams14,aps16,hanbook},
results in Sections~3 and~4 on homogeneous framelets are much more general and complete.
In the last Section~5, we discuss the relations between homogeneous wavelets and nonhomogeneous wavelets with the refinable structure.

\section{Some Auxiliary Results}

In this section we shall recall some basic definitions and results on shift-invariant subspaces of $\dLp{2}$. Then we shall provide some auxiliary results
on Schur decomposition and singular value decomposition for matrices of measurable functions
for studying homogeneous framelets and wavelets with the refinable structure.

A closed subspace $S$ of $\dLp{2}$ is called \emph{shift-invariant} if $f(\cdot-k)\in S$ for all $f\in S$ and $k\in \dZ$.
For a subset $\Phi\subset \dLp{2}$, the shift-invariant space generated by $\Phi$ in $\dLp{2}$ is defined to be
\be \label{sis}
\si(\Phi)=\ol{\mbox{span}\{\phi(\cdot-k) \setsp k\in \dZ, \phi\in \Phi\}}^{\|\cdot\|_{\dLp{2}}},
\ee
where the overhead bar refers to the closure operation in $\dLp{2}$. That is, $\si(\Phi)$ is the smallest closed subspace of $\dLp{2}$ containing the generator set $\Phi$.
Because the space $\dLp{2}$ is separable, every shift-invariant subspace $S$ in $\dLp{2}$ is generated by a countable subset $\Phi\subset \dLp{2}$. By $\#\Phi$ we denote the cardinality of a set $\Phi$. The \emph{length} of a shift-invariant subspace $S$ in $\dLp{2}$ is defined to be
\be \label{len}
\len(S):=\min\{ \#\Phi \setsp S=\si(\Phi)\}
\ee
and the \emph{dimension function} of a shift-invariant space $S$ on $\dR$ is defined to be
\be \label{dimfunc}
\dim_S(\xi):=\dim(\mbox{span}\{ \{\wh{\phi}(\xi+2\pi k)\}_{k\in \dZ} \setsp \phi\in \Phi\}),\qquad \xi\in \dR,
\ee
where $\Phi$ is a countable subset of $\dLp{2}$ satisfying $S=\si(\Phi)$. Up to a set of measure zero, the dimension function is independent of the choice of a generator set $\Phi$ and is a $2\pi\dZ$-periodic measurable function on $\dR$. Moreover, the identity $\len(S)=\|\dim_S(\cdot)\|_{\dTLp{\infty}}$ holds. For a comprehensive study on shift-invariant subspaces in $\dLp{2}$, see de Boor, DeVore and Ron~\cite{bdr94} and Bownik~\cite{bow00} and references therein.

To study shift-invariant subspaces of $\dLp{2}$, we recall the bracket product.
For $f\in (\dLp{2})^{r\times n}$ and $g\in (\dLp{2})^{s\times n}$, we define an $r\times s$ matrix $[f,g]$ of functions by
\be \label{bracket}
[f,g](\xi):=\sum_{k\in \dZ} f(\xi+2\pi k)\ol{g(\xi+2\pi k)}^\tp,\qquad \xi\in \dR
\ee
and an $r\times s$ matrix $\la f,g\ra$ of complex numbers by $\la f,g\ra:=\int_{\dR} f(x)\ol{g(x)}^\tp dx$.
For $f,g\in \dLp{2}$,
it is trivial to check that $\|[\wh{f},\wh{f}]\|_{\dTLp{1}}:=
\frac{1}{(2\pi)^d}\int_{(-\pi,\pi]^d}
[\wh{f},\wh{f}](\xi) d\xi=\frac{1}{(2\pi)^d} \int_{\dR} |\wh{f}(\xi)|^2 d\xi=\|f\|_{\dLp{2}}^2$. Therefore, by the Cauchy-Schwarz inequality, we see that $[\wh{f},\wh{g}](\xi)$ is well defined for almost every $\xi\in \dR$ and $[\wh{f},\wh{g}]\in \dTLp{1}$ by $|[\wh{f},\wh{g}]|^2\le [\wh{f},\wh{f}][\wh{g},\wh{g}]$.
Moreover, one can directly check that the Fourier series of $[\wh{f},\wh{g}]\in \dTLp{1}$ is $\sum_{k\in \dZ} \la f, g(\cdot-k) \ra e^{-ik\cdot \xi}$. Therefore, $\la f, g(\cdot-k)\ra=0$ for all $k\in \dZ\bs\{0\}$ if and only if $[\wh{f},\wh{g}](\xi)=\la f,g\ra$ for almost every $\xi\in \dR$. If in addition $[\wh{f},\wh{g}]\in \dTLp{2}$, then the Parseval's identity yields
\be \label{parseval}
\sum_{k\in \dZ} |\la f, g(\cdot-k)\ra|^2=\frac{1}{(2\pi)^d} \int_{(-\pi,\pi]^d} |[\wh{f},\wh{g}](\xi)|^2 d\xi,\qquad f,g\in \dLp{2}.
\ee
If we only have $[\wh{f},\wh{g}]\in \dTLp{1}$ but $[\wh{f},\wh{g}]\not \in \dTLp{2}$, then the above identity \eqref{parseval} still holds but its both sides are infinity.

Let $\Phi=\{\phi^1,\ldots,\phi^L\}$ be a countable subset of $\dLp{2}$ with $L\in \N \cup\{\infty\}$ ($\Phi$ is countably infinite if $L=\infty$). Throughout this section we shall use the convention $\frac{0}{0}:=0$. Employing the following standard orthogonalization procedure:
\[
\wh{\varphi^\ell}(\xi):=
\begin{cases}
\frac{\wh{\mathring{\phi}^\ell}(\xi)}{
\sqrt{[\wh{\mathring{\phi}^\ell},\wh{\mathring{\phi}^\ell}](\xi)}},
&\text{if $0<[\wh{\mathring{\phi}^\ell},\wh{\mathring{\phi}^\ell}](\xi)<\infty$,}\\
0, &\text{otherwise},
\end{cases}
\]
where
\[
\wh{\mathring{\phi}^\ell}(\xi):=\wh{\phi^\ell}(\xi)-
\sum_{j=1}^{\ell-1} [\wh{\phi^\ell},\wh{\varphi^j}](\xi)\wh{\varphi^j}(\xi),\qquad \ell=1,\ldots,L,
\]
we have $\si(\Phi)=\si(\varphi^1)\oplus\cdots\oplus \si(\varphi^L)$ with $\oplus$ standing for the orthogonal decomposition in $\dLp{2}$
and $[\wh{\varphi^j},\wh{\varphi^j}]=\chi_{\supp([\wh{\varphi^j},\wh{\varphi^j}])}$
for all $j=1,\ldots,L$, where $\supp([\wh{\varphi^j},\wh{\varphi^j}]):=\{\xi\in \dR \setsp [\wh{\varphi^j},\wh{\varphi^j}](\xi)\ne 0\}$.
Cutting $\dR$ by the $2\pi \dZ$-periodic measurable sets $\supp([\wh{\varphi^j},\wh{\varphi^j}]), j=1,\ldots,,L$ and using a simple cut-and-paste technique, we have the following known result (e.g., see \cite{bow00,bdr94} and \cite{hanbook}):

\begin{prop}\label{prop:sis}
Let $\Phi$ be a countable subset of $\dLp{2}$ and define $r:=\len(\si(\Phi))$. Then there exist $\varphi^1,\ldots,\varphi^r\in \si(\Phi)$ such that
\begin{enumerate}
\item[(1)] $\si(\Phi)=\si(\varphi^1)\oplus\cdots\oplus \si(\varphi^r)$ and $\|\varphi^j\|_{\dLp{2}}\le 1$ for all $j=1,\ldots,r$.
\item[(2)]  $[\wh{\varphi^j},\wh{\varphi^j}]=\chi_{\supp([\wh{\varphi^j},\wh{\varphi^j}])}$ holds (i.e., $\{\varphi^j(\cdot-k) \setsp k\in \dZ\}$ is a tight frame of $\si(\varphi^j)$), and $[\wh{\varphi^j},\wh{\varphi^k}]=0$ (i.e., $\si(\varphi^j) \perp \si(\varphi^k)$) for all $j,k=1,\ldots,r$ with $j\ne k$.
\item[(3)] $\dim_{\si(\Phi)}(\xi)=\sum_{j=1}^r [\wh{\varphi^j},\wh{\varphi^j}](\xi)$ for almost every $\xi\in \dR$.
\item[(4)] Up to a set of measure zero, $\supp([\wh{\varphi^{j+1}},\wh{\varphi^{j+1}}])
    \subseteq \supp([\wh{\varphi^j},\wh{\varphi^j}])$ for all $j=1,\ldots,r-1$.
\end{enumerate}
Moreover, for every $f\in \si(\Phi)$, we have
\[
\wh{f}(\xi)=[\wh{f},\wh{\varphi^1}](\xi)\wh{\varphi^1}(\xi)+
\cdots+[\wh{f},\wh{\varphi^r}](\xi)\wh{\varphi^r}(\xi), \qquad a.e.\, \xi\in \dR
\]
and all $[\wh{f},\wh{\varphi^1}],\ldots,[\wh{f},\wh{\varphi^r}]$ are $2\pi\dZ$-periodic measurable functions in $\dTLp{2}$.
\end{prop}

Next, we discuss Hermite matrices of measurable functions.
In this section, the norm $\mnorm{A}$ for a general matrix $A$ stands for the matrix operator norm, that is, $\mnorm{A}:=\sup_{\|x\|_{l_2}\le 1} \|Ax\|_{l_2}$.
An $r\times r$ square matrix $A$ is called \emph{a Hermite matrix} if $A=\ol{A}^\tp$. We now study Hermite matrices of measurable functions. To do so, we need the following result.

\begin{lemma}\label{lem:kernel}
Let $E$ be a (Lebesgue) measurable subset of $\dR$ and $A: E\rightarrow \C^{r\times s}$ be an $r\times s$ matrix of measurable functions on $E$. If the rank of $A(\xi)$ is a constant integer $n$ for all $\xi\in E$, then there exists an $s\times (s-n)$ matrix $V: E\rightarrow \C^{s\times (s-n)}$ of measurable functions on $E$ such that
\be \label{V:kernel}
A(\xi) V(\xi)=0\quad \mbox{and}\quad
\ol{V(\xi)}^\tp V(\xi)=I_{s-n},\qquad \forall\; \xi\in E.
\ee
\end{lemma}

\bp We write the $r\times s$ matrix $A(\xi)$ as follows:
\[
A(\xi)=\left[ \begin{matrix}A_1(\xi) &A_2(\xi)\\
A_3(\xi) &A_4(\xi)\end{matrix}\right],\qquad \xi\in E,
\]
where $A_1$ is an $n\times n$ matrix.
A matrix has rank $n$ if and only if all its $(n+1)\times(n+1)$ minors are zero and at least one of its $n\times n$ minors is nonzero.
Therefore, the measurable set $E$ can be written as a disjoint union of finitely many measurable subsets $F_1,\ldots,F_m$ such that on each $F\in \{F_1,\ldots,F_m\}$, there is a fixed $n\times n$ minor which does not vanish on $F$. Permutating rows and columns of $A(\xi)$, without loss of generality, we can assume that $\det(A_1(\xi))\ne 0$ for all $\xi\in F$. Define an $s\times(s-n)$ matrix $V$ by
\[
V(\xi):=\left[ \begin{matrix} -K(\xi)\\ I_{s-n}\end{matrix}\right]M^{-1/2}(\xi)
\quad\mbox{with}\quad
K(\xi):=A_1^{-1}(\xi)A_2(\xi),
M(\xi):=I_{s-n}+\ol{K(\xi)}^\tp K(\xi),
\]
for $\xi\in F$.
We now prove that $V$ is a well-defined desired matrix on $F$. It is trivial to check that
\[
\ol{V(\xi)}^\tp V(\xi)=
M^{-1/2}(\xi)(I_{s-n}+\ol{K(\xi)}^\tp K(\xi))M^{-1/2}(\xi)=M^{-1/2}(\xi)M(\xi) M^{-1/2}(\xi)
=I_{s-n}.
\]
To prove $A(\xi) V(\xi)=0$ for all $\xi\in F$, we observe that
\[
\left[
\begin{matrix} I_n &0\\ -A_3(\xi) A_1^{-1}(\xi) &I_{r-n}\end{matrix}\right]
\left[ \begin{matrix}A_1(\xi) &A_2(\xi)\\
A_3(\xi) &A_4(\xi)\end{matrix}\right]
=\left[ \begin{matrix}A_1(\xi) &A_2(\xi)\\
0 &A_4(\xi)-A_3(\xi)A_1^{-1}(\xi)A_2(\xi)
\end{matrix}\right].
\]
Since the matrix $A(\xi)$ has rank $n$ and the $n\times n$ minor $\det(A_1(\xi))\ne 0$ for all $\xi\in F$, we must have $A_4(\xi)-A_3(\xi)A_1^{-1}(\xi)A_2(\xi)=0$ for all $\xi\in F$.
By the definition of the matrix $V(\xi)$, it is now straightforward to directly check that $A(\xi)V(\xi)=0$ for all $\xi\in F$.

To prove that all the entries of the matrix $V$ are measurable, since $\det(M(\xi)) \ge 1$,
it suffices to prove that $M^{1/2}$ is measurable.
Approximating the entries of $K$ by measurable simple functions, we see that there is a sequence $\{M_j\}_{j\in \N}$ of $(s-n)\times (s-n)$ Hermite matrices of measurable simple functions on $F$ such that $M_j(\xi)\ge I_{s-n}$ and $\lim_{j\to \infty} \mnorm{M_j(\xi)-M(\xi)}=0$ for all $\xi\in F$.
Define $N_j(\xi):=M_j^{1/2}(\xi)-M^{1/2}(\xi)$.
Let $\gl_{\max}(\xi)$ be the largest eigenvalue of $N_j(\xi)$ in modulus having an eigenvector $v(\xi)\in \C^r$ satisfying $\|v(\xi)\|_{l_2}=1$ and $N_j(\xi)v(\xi)=\gl_{\max}(\xi)v(\xi)$.
Since $N_j(\xi)$ is a Hermite matrix, its eigenvalue $\gl_{\max}(\xi)$ must be real-valued.
Then $\mnorm{N_j(\xi)}=|\gl_{\max}(\xi)|$
and
\begin{align*}
\mnorm{M_j(\xi)-M(\xi)}
&=\mnorm{M_j^{1/2}(\xi)N_j(\xi)+N_j(\xi) M^{1/2}(\xi)} \\
&\ge \left|\ol{v(\xi)}^\tp M_j^{1/2}(\xi) N_j(\xi)v(\xi)+
\ol{v(\xi)}^\tp N_j(\xi) M^{1/2}(\xi) v(\xi)\right|\\
&=\left|\gl_{\max}(\xi)\right|\left
(\ol{v(\xi)}^\tp M_j^{1/2}(\xi)v(\xi)+
\ol{v(\xi)}^\tp M^{1/2}(\xi)v(\xi)\right)\\
&\ge 2|\gl_{\max}(\xi)|=2\mnorm{N_j(\xi)},
\end{align*}
where we used $M_j\ge I_{s-n}$, $M\ge I_{s-n}$, and $\|v(\xi)\|_{l_2}=1$. Consequently, we have $\mnorm{M_j^{1/2}(\xi)-M^{1/2}(\xi)}
=\mnorm{N_j(\xi)}\le \frac{1}{2} \mnorm{M_j(\xi)-M(\xi)}\to 0$ as $j\to \infty$.
That is, we proved $\lim_{j\to \infty} \mnorm{M_j^{1/2}(\xi)-M^{1/2}(\xi)}=0$ for all $\xi\in F$.
Since $M_j$ is a Hermite matrix of measurable simple functions, all the entries of the matrix $M_j^{1/2}$ are measurable for every $j\in \N$. Consequently, all the entries of the matrix $M^{1/2}$ are measurable by $\lim_{j\to \infty}M_j^{1/2}(\xi)=M^{1/2}(\xi)$ for all $\xi\in F$. Hence, all the entries of the matrix $V$ are measurable.
The proof is completed by putting all the pieces of $V$ together on $F_1,\ldots,F_m$.
\ep

For an $r\times r$ Hermite matrix $A$,
the well-known Courant-Fisher Theorem tells us that
\[
\gl_j(A)=\max_{\dim(V)=j}\min_{x\in V\bs\{0\}} \frac{\ol{x}^\tp A x}{\|x\|^2_{l_2}},
\]
where $\gl_j(A)$ is the $j$th largest eigenvalue of $A$. Let $B$ be another $r\times r$ Hermite matrix.
Note that
\[
|\ol{x}^\tp A x-\ol{x}^\tp Bx|= |\ol{x}^\tp(A-B)x|\le \mnorm{A-B}\|x\|^2_{l_2}.
\]
Now we can deduce trivially from the Courant-Fisher Theorem that
\be \label{eig:AB}
|\gl_j(A)-\gl_j(B)|\le \mnorm{A-B}.
\ee

To study homogeneous framelets and wavelets with the refinable structure, we need the following result on symmetric Schur decomposition for Hermite matrices of measurable functions, for which we shall provide a simple proof here.

\begin{theorem}\label{thm:schur}
Let $A: \dR \rightarrow \C^{r\times r}$ be an $r\times r$ matrix of (Lebesgue) measurable functions on $\dR$ such that $\ol{A(\xi)}^\tp=A(\xi)$ for all $\xi\in \dR$. Let $\gl_j(\xi)$ be the $j$th largest eigenvalue of $A(\xi)$ for $j=1,\ldots,r$ ordered nonincreasingly by $\gl_1(\xi)\ge \gl_2(\xi)\ge \cdots\ge \gl_r(\xi)$.
Then all the eigenvalue functions $\gl_1,\ldots,\gl_r$ are measurable and there exists an $r\times r$ unitary matrix $U$ of measurable functions on $\dR$ such that
\be \label{evd}
A(\xi)=U(\xi)\diag(\gl_1(\xi),\ldots,\gl_r(\xi)) \ol{U(\xi)}^\tp \quad \mbox{and}\quad
\ol{U(\xi)}^\tp U(\xi)=I_r,\qquad \forall\, \xi \in \dR.
\ee
\end{theorem}

\bp Since $\ol{A(\xi)}^\tp=A(\xi)$, the matrix $A(\xi)$ is a Hermite matrix and it is well known that all its eigenvalues $\gl_1(\xi),\ldots,\gl_r(\xi)$ are real numbers. We first prove that $\gl_1,\ldots,\gl_r$ are measurable.
Let $\{A_n(\xi)\}_{n\in \N}$ be a sequence of $r\times r$ matrices of measurable simple functions on $\dR$ such that
$\ol{A_n(\xi)}^\tp=A_n(\xi)$ and
$\lim_{n\to \infty} A_n(\xi)=A(\xi)$ for every $\xi\in \dR$.
Obviously, for every $j=1,\ldots,r$, the $j$th largest eigenvalue $\gl_j(A_n(\xi))$ of $A_n(\xi)$ is a measurable simple function. Now it follows directly from the inequality \eqref{eig:AB} that $\gl_j(A(\xi))=\lim_{n\to\infty} \gl_j(A_n(\xi))$. Hence, all the functions $\gl_j(A(\xi))$ (i.e., $\gl_j(\xi)$) must be measurable.

We now prove the existence of a measurable unitary matrix $U$. Since all $\gl_j$ are measurable, we can write $\dR$ as a disjoint union of finitely many measurable sets $E_1,\ldots,E_m$ such that each distinct eigenvalue $\gl_j(\xi)$ has the same multiplicity for all $\xi$ within each  $E\in \{E_1,\ldots, E_m\}$, that is,
\be \label{gl:partition}
\gl_1(\xi)=\cdots=\gl_{k_1}(\xi)>\gl_{k_1+1}(\xi)=\cdots=
\gl_{k_2}(\xi)>\gl_{k_2+1}(\xi)=\cdots > \gl_{k_n+1}(\xi)=\cdots=\gl_r(\xi)
\ee
so that the integers $k_1,\ldots,k_n$ are independent of $\xi\in E$.
Now we consider each of the distinct eigenvalues $\gl_{k_1},\gl_{k_2},\ldots,\gl_{k_n},\gl_r$.
For simplicity, it suffices to consider $\gl_{k_1}$, which is just $\gl_1$. By our assumption, the matrix $\gl_1(\xi)I_r-A(\xi)$ has rank $r-k_1$ for all $\xi\in E$.
By Lemma~\ref{lem:kernel}, there exists an $r\times k_1$ matrix $V_1$ of measurable functions on $E$ such that
$\ol{V_1(\xi)}^\tp V_1(\xi)=I_{k_1}$ and
$[\gl_1(\xi)I_r-A(\xi)] V_1(\xi)=0$ (i.e., $A(\xi)V_1(\xi)=\gl_1(\xi) V_1(\xi)$) for all $\xi\in E$.
Performing this procedure for each distinct eigenvalue of $\gl_{k_1},\gl_{k_2},\ldots,\gl_{k_n},\gl_r$ and putting all the $V$'s together as an $r\times r$ matrix $U$, we conclude that $A(\xi) U(\xi)=U(\xi)\mbox{diag}(\gl_1(\xi),\ldots,\gl_r(\xi))$
and $\ol{U(\xi)}^\tp U(\xi)=I_r$ for all $\xi\in E$, where we used the fact that $\ol{V_j(\xi)}^\tp V_\ell(\xi)=0$ for two different eigenvalues $\gl_{k_j}(\xi)\ne \gl_{k_\ell}(\xi)$.
The proof is completed by putting together all the pieces of $U$ on the measurable subsets $E_1,\ldots,E_m$.
\ep

As a direct consequence of Theorem~\ref{thm:schur} and Lemma~\ref{lem:kernel}, we also have the following result on singular value decomposition of a matrix of measurable functions.

\begin{cor}\label{cor:svd}
Let $A: \dR\rightarrow \C^{r\times s}$ be an $r\times s$ matrix of measurable functions on $\dR$. Let $\sigma_j(\xi)$ be the $j$th largest singular value of $A(\xi)$ ordered nonincreasingly by $\sigma_1(\xi)\ge \sigma_2(\xi)\ge \cdots\ge \sigma_{\min(r,s)}(\xi)\ge 0$. Then all $\sigma_1,\ldots,\sigma_{\min(r,s)}$ are nonnegative measurable functions and there exist an $r\times r$ unitary matrix $U$ and an $s\times s$ unitary matrix $V$ of measurable functions on $\dR$ such that $\ol{U(\xi)}^\tp U(\xi)=I_r$, $\ol{V(\xi)}^\tp V(\xi)=I_s$,
and the $r\times s$ matrix
$\ol{U(\xi)}^\tp A(\xi) V(\xi)$
is an $r\times s$ diagonal matrix with the first $\min(r,s)$ diagonal entries being
$\sigma_1(\xi),\ldots,\sigma_{\min(r,s)}(\xi)$ and with all its other entries being zero.
\end{cor}

\bp Without loss of generality, we assume $r\le s$. Note that $\sigma_1^2,\ldots,\sigma_r^2$ are all the eigenvalues of the Hermite matrix $A(\xi)\ol{A(\xi)}^\tp$.
By Theorem~\ref{thm:schur}, all the functions $\sigma_1^2,\ldots,\sigma_r^2$ are measurable and there exists an $r\times r$ unitary matrix $U$ of measurable functions such that $\ol{U(\xi)}^\tp U(\xi)=I_r$ and
$\ol{U(\xi)}^\tp A(\xi)\ol{A(\xi)}^\tp U(\xi)=\mbox{diag}(\sigma_1^2(\xi),\ldots,\sigma_r^2(\xi))$.
Consequently, all $\sigma_1,\ldots,\sigma_r$ are nonnegative measurable functions.
We define $D$ to be the $r\times s$ diagonal matrix whose first $r$ diagonal entries are
$\sigma_1^{-1}(\xi),\ldots,\sigma^{-1}_r(\xi)$, under the convention that $\sigma_j^{-1}(\xi):=0$ if $\sigma_j(\xi)=0$, and all the other entries of the matrix $D$ are zero.
Define an $s\times s$ matrix $\tilde{V}(\xi):=\ol{A(\xi)}^\tp U(\xi)D(\xi)$ of measurable functions.
By $A(\xi)\ol{A(\xi)}^\tp=
U(\xi)\mbox{diag}(\sigma_1^2(\xi),\ldots,\sigma_r^2(\xi))\ol{U(\xi)}^\tp$, since $U$ is unitary, we have
\be \label{Vsquare}
\begin{split}
\ol{\tilde{V}(\xi)}^\tp \tilde{V}(\xi)&=
(D(\xi))^\tp\ol{U(\xi)}^\tp A(\xi)\ol{A(\xi)}^\tp U(\xi) D(\xi)
=(D(\xi))^\tp \mbox{diag}(\sigma_1^2(\xi),\ldots,\sigma_r^2(\xi))
D(\xi)\\
&=\mbox{diag}(\mbox{sgn}(\sigma_1(\xi)),\ldots,
\mbox{sgn}(\sigma_r(\xi)),0,\ldots,0),
\end{split}
\ee
where
\be \label{sgn}
\mbox{sgn}(c):=1,\quad \mbox{if}\; c>0, \qquad \mbox{sgn}(0):=0, \qquad \mbox{and} \quad \mbox{sgn}(c)=-1,\quad \mbox{if}\; c<0.
\ee
This shows that the columns of $\tilde{V}$ are mutually orthogonal and the $l_2$-norm of each column of $\tilde{V}$ is either $1$ or $0$.
Therefore, the set $\dR$ can be written as a disjoint union of finitely many measurable subsets $E_1,\ldots,E_m$ such that on each $E\in \{E_1,\ldots,E_m\}$, exactly $n$ fixed columns of $\tilde{V}(\xi)$ have $l_2$-norm one (i.e., the integer $n$ is independent of $\xi\in E$). Without loss of generality, we can assume that $\tilde{V}(\xi)=[V_1(\xi),0]$ such that the $s\times n$ matrix $V_1$ satisfies $\ol{V_1(\xi)}^\tp V_1(\xi)=I_n$ for all $\xi\in E$ (i.e., every column of $V_1(\xi)$ has $l_2$-norm one).
Now it follows directly from \eqref{Vsquare} that
$\sigma_1(\xi)\ge \cdots\ge \sigma_n(\xi)>0$ and $\sigma_{n+1}(\xi)=\cdots=\sigma_r(\xi)=0$
for all $\xi\in E$. Moreover, by the definition of $\tilde{V}(\xi)=\ol{A(\xi)}^\tp U(\xi) D(\xi)$, we see that
the $r\times s$ matrix
\begin{align*}
\ol{U(\xi)}^\tp A(\xi)\tilde{V}(\xi)
&=
\ol{U(\xi)}^\tp A(\xi)\ol{A(\xi)}^\tp U(\xi) D(\xi)=\mbox{diag}(
\sigma_1^2(\xi),\ldots,\sigma_r^2(\xi))D(\xi)\\
&=[\mbox{diag}(\sigma_1(\xi),\ldots,\sigma_n(\xi),0,\ldots,0),0]=
[\mbox{diag}(\sigma_{1}(\xi),\ldots,\sigma_r(\xi)),0].
\end{align*}
Since $\ol{U(\xi)}^\tp A(\xi)\ol{A(\xi)}^\tp U(\xi)=\mbox{diag}(\sigma^2_1(\xi),\ldots,\sigma^2_n(\xi),0,\ldots,0)$ with $\sigma_1(\xi)\ge \cdots\ge \sigma_n(\xi)>0$ for all $\xi\in E$, the rank of the $r\times s$ matrix $\ol{U(\xi)}^\tp A(\xi)$ is $n$ for all $\xi\in E$.
By Lemma~\ref{lem:kernel}, there exists an $s\times (s-n)$ matrix $V_2(\xi)$ of measurable functions on $E$ such that $\ol{U(\xi)}^\tp A(\xi) V_2(\xi)=0$ and $\ol{V_2(\xi)}^\tp V_2(\xi)=I_{s-n}$ for all $\xi \in E$. Define an $s\times s$ matrix $V(\xi):=[V_1(\xi),V_2(\xi)]$. Since $\ol{U(\xi)}^\tp A(\xi) V_2(\xi)=0$, we trivially have
\[
\ol{U(\xi)}^\tp A(\xi)V(\xi)
=\ol{U(\xi)}^\tp A(\xi)\tilde{V}(\xi)=
[\mbox{diag}(\sigma_{1}(\xi),\ldots,\sigma_r(\xi)),0],
\]
for all $\xi\in E$.
By $\tilde{V}(\xi)=[V_1(\xi),0]=\ol{A(\xi)}^\tp U(\xi) D(\xi)$ and $\ol{U(\xi)}^\tp A(\xi)V_2(\xi)=0$,
we have
\[
\ol{[V_1(\xi), 0]}^\tp V_2(\xi)=\ol{\tilde{V}(\xi)}^\tp V_2(\xi)=
(D(\xi))^\tp\, \ol{U(\xi)}^\tp A(\xi) V_2(\xi)=0.
\]
Therefore, we must have $\ol{V_1(\xi)}^\tp V_2(\xi)=0$ and hence $\ol{V(\xi)}^\tp V(\xi)=I_s$ for all $\xi\in E$.
The proof is completed by putting $V$ together on all the pieces of measurable sets $E_1,\ldots,E_m$.
\ep

\section{Link Homogeneous Framelets to Nonhomogeneous Framelets}

In this section we study homogeneous framelets and homogeneous tight framelets by linking them to nonhomogeneous framelets and tight framelets with the refinable structure.

To study homogeneous framelets, we need the following result.

\begin{lemma}\label{lem:hf}
Let $H$ be a countable subset of $\dLp{2}$. Assume that $\{h(\cdot-k) \setsp k\in \dZ, h\in H\}$ is a Bessel sequence in $\dLp{2}$, i.e.,
there exists a positive constant $C$ such that
\be \label{bessel}
\sum_{h\in H} \sum_{k\in \dZ} |\la f, h(\cdot-k)\ra|^2\le C\|f\|_{\dLp{2}}^2, \qquad \forall\, f\in \dLp{2}.
\ee
If $r:=\len(\si(H))<\infty$, then there exists $\eta^1,\ldots,\eta^r\in \si(H)$ such that
\be \label{basic:identity}
\sum_{\ell=1}^r \sum_{k\in \dZ} |\la f, \eta^\ell(\cdot-k)\ra|^2=
\sum_{h \in H} \sum_{k\in \dZ} |\la f, h(\cdot-k)\ra|^2,\qquad \forall\, f\in \dLp{2}
\ee
and $\si(\{\eta^1,\ldots,\eta^r\})=\si(H)$.
\end{lemma}

\bp Since $r=\len(\si(H))$, by Proposition~\ref{prop:sis},
there exist $\varphi^1,\ldots,\varphi^r\in \si(H)$ such that items (1)--(4) in Proposition~\ref{prop:sis} hold with $\Phi$ being replaced by $H$.
Therefore, for each $h\in H\subset \si(H)$, we have
\[
\wh{h}=u_{h,1}\wh{\varphi^1}+\cdots+u_{h,r} \wh{\varphi^r}\quad \mbox{with}\quad
u_{h,j}:=[\wh{h},\wh{\varphi^j}]\in \dTLp{2},\qquad j=1,\ldots,r.
\]
%
%
For $f\in \dLp{2}$, by the identity \eqref{parseval} we have
\begin{align*}
(2\pi)^d \sum_{h\in H} \sum_{k\in \dZ}
|\la f, h(\cdot-k)\ra|^2&=
\int_{(-\pi,\pi]^d} \sum_{h\in H}
|[\wh{f},\wh{h}](\xi)|^2 d\xi\\
&=\int_{(-\pi,\pi]^d} \sum_{h\in H}
|[\wh{f},u_{h,1}\wh{\varphi^1}+\cdots+u_{h,r}\wh{\varphi^r}](\xi)|^2 d\xi\\
&=\int_{(-\pi,\pi]^d} \sum_{h\in H}
\sum_{j=1}^r \sum_{k=1}^r
[\wh{f},\wh{\varphi^j}](\xi)
\ol{u_{h,j}(\xi)}u_{h,k}(\xi)[\wh{\varphi^k},\wh{f}](\xi)
d\xi.
\end{align*}
If $f=\varphi^\ell$, then we deduce from the above identity and the fact $[\wh{\varphi^j},\wh{\varphi^k}]=0$ for all $j\ne k$ that
\[
\int_{(-\pi,\pi]^d}
\sum_{h\in H} |u_{h,\ell}(\xi)|^2 d\xi=
\int_{(-\pi,\pi]^d}
\sum_{h\in H} |[\wh{\varphi^\ell},\wh{\varphi^\ell}](\xi)|^2
|u_{h,\ell}(\xi)|^2 d\xi
=(2\pi)^d\sum_{h\in H}\sum_{k\in \dZ} |\la \varphi^\ell,h(\cdot-k)\ra|^2,
\]
where we used the facts that $[\wh{\varphi^\ell},\wh{\varphi^\ell}]=
\chi_{\supp([\wh{\varphi^\ell},\wh{\varphi^\ell}])}$
and $u_{h,\ell}(\xi)=[\wh{h},\wh{\varphi^\ell}](\xi)=0$
for $\xi\not \in \supp([\wh{\varphi^\ell},\wh{\varphi^\ell}])$.
Consequently, by our assumption in \eqref{bessel} and $\|\varphi^\ell\|_{\dLp{2}}\le 1$, we proved
$\sum_{h\in H} |u_{h,\ell}|^2\in \dTLp{1}$ for all $\ell=1,\ldots,r$. Therefore, the following $r\times r$ matrix $A(\xi)$ with its $(j,k)$-entry
\be \label{matrix:A}
[A(\xi)]_{j,k}:=\sum_{h\in H} \ol{u_{h,j}(\xi)}u_{h,k}(\xi),\qquad \xi\in \dR, j,k=1,\ldots,r
\ee
is well defined with all entries belonging to $\dTLp{1}$. Since $\ol{A(\xi)}^\tp=A(\xi)$, by Theorem~\ref{thm:schur}, there exists an $r\times r$ unitary matrix $U$ of $2\pi\dZ$-periodic measurable functions such that \eqref{evd} holds with $\gl_j(\xi)$ being the $j$th largest eigenvalue of $A(\xi)$.
Since $A(\xi)\ge 0$ and each entry of $A$ belongs to $\dTLp{1}$, we must have $0\le \gl_r\le \cdots\le \gl_2 \le \gl_1$ and $\sqrt{\gl_1}\in \dTLp{2}$ by \eqref{evd}.
Define $A^{1/2}(\xi)=U(\xi)\diag(\sqrt{\gl_1(\xi)},
\ldots,\sqrt{\gl_r(\xi)})
\ol{U(\xi)}^\tp$ and
\be \label{def:eta:A}
\wh{\eta}(\xi):=
A^{1/2}(\xi)\wh{\varphi}(\xi)
\quad \mbox{with}
\quad
\eta:=(\eta^1,\ldots,\eta^r)^\tp,
\varphi:=(\varphi^1,\ldots,\varphi^r)^\tp.
\ee
Then by item (2) of Proposition~\ref{prop:sis}, we have $[\wh{\varphi},\wh{\varphi}](\xi)\le I_r$ and
\[
[\wh{\eta},\wh{\eta}](\xi)=
A^{1/2}(\xi)[\wh{\varphi},\wh{\varphi}](\xi) A^{1/2}(\xi)
\le A^{1/2}(\xi)A^{1/2}(\xi)=A(\xi).
\]
Therefore, $\sum_{j=1}^r [\wh{\eta^j},\wh{\eta^j}](\xi)
=\mbox{trace}([\wh{\eta},\wh{\eta}](\xi))
\le \mbox{trace}(A(\xi))$.
By $A\in (\dTLp{1})^{r\times r}$, we conclude that  $ \sum_{j=1}^r [\wh{\eta^j},\wh{\eta^j}] \in \dTLp{1}$ and consequently,
$\eta^1,\ldots,\eta^r\in \dLp{2}$ by $\|\eta^j\|_{\dLp{2}}^2=\|[\wh{\eta^j},\wh{\eta^j}]\|_{\dTLp{1}}$.
Moreover, we have
\begin{align*}
(2\pi)^d \sum_{h\in H} \sum_{k\in \dZ}
|\la f, h(\cdot-k)\ra|^2
&=\int_{(-\pi,\pi]^d}
\sum_{j=1}^r \sum_{k=1}^r
[\wh{f},\wh{\varphi^j}]
[A(\xi)]_{j,k}
[\wh{\varphi^k},\wh{f}](\xi)
d\xi\\
&=\int_{(-\pi,\pi]^d}
[\wh{f},\wh{\varphi}](\xi)
A^{1/2}(\xi) A^{1/2}(\xi)
[\wh{\varphi},\wh{f}](\xi)
d\xi\\
&=
\int_{(-\pi,\pi]^d}
[\wh{f},\wh{\eta}](\xi)
[\wh{\eta},\wh{f}](\xi)
d\xi\\
&=
\int_{(-\pi,\pi]^d}
\sum_{j=1}^r |[\wh{f},\wh{\eta^j}](\xi)|^2
d\xi
=(2\pi)^d \sum_{j=1}^r \sum_{k\in \dZ} |\la f,\eta^j(\cdot-k)\ra|^2.
\end{align*}
This proves \eqref{basic:identity}.
Now it follows directly from \eqref{basic:identity} that $f \perp \si(H)$ if and only if $f \perp \si(\{\eta^1,\ldots,\eta^r\})$. This proves $\si(\{\eta^1,\ldots,\eta^r\})=\si(H)$.
\ep

As a direct consequence of Lemma~\ref{lem:hf}, we have

\begin{cor}\label{cor:hf:1}
Let $\dm$ be a $d\times d$ invertible real-valued matrix.
Let $\Psi$ be a countable subset of $\dLp{2}$ such that $s:=\len(\si(\Psi))<\infty$.
If $\Psi$ is a homogeneous $\dm$-framelet in $\dLp{2}$ (i.e., $\AS(\Psi)$ is a frame for $\dLp{2}$ satisfying \eqref{hframelet}),
then there exists a subset $H:=\{\eta^1,\ldots,\eta^s\}\subset \si(\Psi)$ such that $\si(H)=\si(\Psi)$, $H$ is a homogeneous $\dm$-framelet in $\dLp{2}$ with the same lower and upper frame bounds, and
\be \label{framelet:preserve}
\sum_{\psi\in \Psi} \sum_{k\in \dZ}|\la f,\psi_{\dm^j;k}\ra|^2
=\sum_{\ell=1}^s \sum_{k\in \dZ}
|\la f, \eta^\ell_{\dm^j;k}\ra|^2,\qquad \forall\, f\in \dLp{2}, j\in \Z.
\ee
\end{cor}

\bp The claim follows directly from Lemma~\ref{lem:hf} and the simple fact $\la f,\psi_{\dm^j;k}\ra=\la f_{\dm^{-j};0},\psi(\cdot-k)\ra$ for all $j\in \Z$ and $k\in \dZ$.
\ep

We now connect a homogeneous framelet with a nonhomogeneous framelet as follows.

\begin{prop}\label{prop:hf:1}
Let $\dm$ be a $d\times d$ expansive integer matrix.
Let $\Psi$ be a countable subset of $\dLp{2}$ such that
$\sum_{\psi\in \Psi} \|\psi\|^2_{\dLp{2}}<\infty$.
Define
\be \label{Psi:H}
H:=\{ |\det(\dm)|^{-j} \psi(\dm^{-j}\cdot) \setsp j\in \N, \psi\in \Psi\}.
\ee
Suppose that $r:=\len(\si(H))<\infty$. Then $\Psi$ is a homogeneous $\dm$-framelet in $\dLp{2}$ satisfying \eqref{hframelet} if and only if there exists a subset $\Phi:=\{\varphi^1,\ldots,\varphi^r\}\subset \si(H)$ such that
$\si(\Phi)=\si(H)$,
\be \label{H:Phi}
\sum_{\varphi\in \Phi} \sum_{k\in \dZ}
|\la f,\varphi(\cdot-k)\ra|^2=
\sum_{j=1}^\infty \sum_{\psi\in \Psi}\sum_{k\in \dZ} |\la f, |\det(\dm)|^{-j} \psi(\dm^{-j}(\cdot-k))\ra|^2,
\quad \forall\; f\in \dLp{2},
\ee
and $\{\Phi; \Psi\}$ is an $\dm$-framelet in $\dLp{2}$ with the same lower and upper frame bounds, i.e., for all $f\in \dLp{2}$,
\be \label{framelet}
C_1\|f\|^2_{\dLp{2}}\le
\sum_{\varphi\in \Phi}
\sum_{k\in \dZ} |\la f, \varphi(\cdot-k)\ra|^2+
\sum_{j=0}^\infty \sum_{\psi\in \Psi}\sum_{k\in \dZ} |\la f,\psi_{\dm^j;k}\ra|^2\le C_2\|f\|_{\dLp{2}}^2.
\ee
\end{prop}

\bp Sufficiency ($\Rightarrow$). By \cite[Theorem~5.5]{rs97} and \cite[Theorem~2]{css98} for a finite set $\Psi$ and \cite[Theorem~4.3.4]{hanbook} for a countable set $\Psi$ satisfying $\sum_{\psi\in \Psi} \|\psi\|^2_{\dLp{2}}<\infty$,
the inequality \eqref{framelet} holds with $\Phi$ being replaced by $H$ and having the same lower and upper frame bounds. Now the claim follows directly from Lemma~\ref{lem:hf}.

Necessity ($\Leftarrow$).
If $\{\Phi;\Psi\}$ is an $\dm$-framelet in $\dLp{2}$, since $\dm$ is expansive,
then it has been proved in \cite[Proposition~4]{han12} that $\Psi$ must be a homogeneous $\dm$-framelet in $\dLp{2}$.
\ep

It is of interest to know when the condition $\len(\si(H))<\infty$ in Proposition~\ref{prop:hf:1} is satisfied
where $H$ is defined in \eqref{Psi:H}. For some special subsets $\Psi$ of $\dLp{2}$, we have

\begin{lemma}
Let $\Psi$ be a finite subset of $\dLp{2}$. If there exists $\gep>0$ such that
$\supp(\wh{\psi})\subseteq K_\gep$ for all $\psi\in \Psi$
with $K_\gep:=\{\xi\in \dR \setsp \gep\le \|\xi\|\le \gep^{-1}\}$, then $\len(\si(H))<\infty$ with $H$ defined in \eqref{Psi:H}.
\end{lemma}

\bp Note that the Fourier transform of $|\det(\dm)|^{-j} \psi(\dm^{-j}\cdot)$ is $\wh{\psi}((\dm^\tp)^j\cdot)$. Therefore, by the definition of the dimension function in \eqref{dimfunc}, we have
\[
\dim_{\si(H)}(\xi)=\dim(\mbox{span}\{\{\wh{\psi}((\dm^\tp)^j(\xi+2\pi k))\}_{k\in \dZ} \setsp j\in \N, \psi\in \Psi\}).
\]
Since each $\wh{\psi}$ vanishes outside $K_\gep$, we have
\be \label{supp:psi}
\supp(\wh{\psi}((\dm^\tp)^j\cdot))\subseteq (\dm^\tp)^{-j}\supp(\wh{\psi})
\subseteq (\dm^\tp)^{-j} K_\gep.
\ee
Because $K_\gep$ is a bounded set and $\dR\bs K_\gep$ contains a neighborhood of the origin, the above inequalities in \eqref{supp:psi} imply that there exists $J\in \N$ such that for all $j\ge J$, $\supp(\wh{\psi}((\dm^\tp)^j\cdot))$ is contained inside $(-\pi/2,\pi/2)^d$ and $\wh{\psi}((\dm^\tp)^j(\xi+2\pi k))=0$ for all $k\in \dZ\bs\{0\}$ and $\xi\in (-\pi,\pi]^d$.
On the other hand, since $K_\gep=\{\xi\in \dR \setsp \gep\le \|\xi\|\le \gep^{-1}\}$ and $\dm$ is expansive, we trivially deduce from \eqref{supp:psi} that
\[
\sum_{j\in \N} \chi_{\supp(\wh{\psi}((\dm^\tp)^j\cdot))}(\xi)
\le \sum_{j\in \N} \chi_{(\dm^\tp)^{-j} K_\gep}(\xi)
=\sum_{j\in \N} \chi_{K_\gep}((\dm^\tp)^j \xi)\in \dLp{\infty}.
\]
Consequently, for $\xi \in (-\pi,\pi]^d$, by $\#\Psi<\infty$ (i.e., $\Psi$ is a finite set), we must have
\begin{align*}
\dim_{\si(H)}(\xi)
&\le
\dim(\mbox{span}\{\{\wh{\psi}((\dm^\tp)^j(\xi+2\pi k))\}_{k\in \dZ} \setsp 1\le j<J, \psi\in \Psi\})\\
&\qquad +
\dim(\mbox{span}\{\{\wh{\psi}((\dm^\tp)^j(\xi+2\pi k))\}_{k\in \dZ} \setsp j\ge J, \psi\in \Psi\})\\
&\le
J(\#\Psi)+
\sum_{\psi\in \Psi}\left\| \sum_{j=J}^\infty \chi_{\supp(\wh{\psi}((\dm^\tp)^j\cdot))}(\cdot)\right\|_{\dLp{\infty}}\\
&\le J(\#\Psi)+(\#\Psi)\left\| \sum_{j\in \N} \chi_{K_\gep}((\dm^\tp)^j \cdot)\right\|_{\dLp{\infty}}<\infty.
\end{align*}
This proves that $\len(\si(H))=\|\dim_{\si(H)}(\cdot)\|_{\dLp{\infty}}<\infty$.
\ep

For a $d\times d$ real-valued matrix $U$ and a subset $\Phi\subset \dLp{2}$, we define
\be \label{sisU}
\si_U(\Phi):=\{f(U\cdot) \setsp f\in \si(\Phi)\}=
\ol{\mbox{span}\{\phi(U\cdot-k) \setsp k\in \dZ, \phi\in \Phi\}}^{\|\cdot\|_{\dLp{2}}}.
\ee
In terms of shift-invariant spaces, the refinable structure in \eqref{refstr} is equivalent to $\si(\Phi)\cup \si(\Psi)\subseteq \si_\dm(\Phi)$, where $\Phi$ and $\Psi$ are subsets containing all the entries in the vector functions $\phi$ and $\psi$, respectively.

For homogeneous framelets having the refinable structure in \eqref{refstr}, we have

\begin{theorem}\label{thm:hf:refstr}
Let $\dm$ be a $d\times d$ expansive integer matrix.
Let $\Phi=\{\phi^1,\ldots,\phi^r\}$ be a finite subset of $\dLp{2}$ and $\Psi$ be a countable subset of $\dLp{2}$ such that $\sum_{\psi\in \Psi} \|\psi\|_{\dLp{2}}^2<\infty$ and the refinable structure $\si(\Phi)\cup\si(\Psi)\subseteq \si_\dm(\Phi)$ holds (i.e., \eqref{refstr} holds with $\phi:=(\phi^1,\ldots,\phi^r)^\tp$
and $\psi$ being the column vector by listing all the elements in $\Psi$
for some matrices $\wh{a}$ and $\wh{b}$ of $2\pi\dZ$-periodic measurable functions).
If $\Psi$ is a homogeneous $\dm$-framelet in $\dLp{2}$ satisfying \eqref{hframelet}, then there exist subsets $\mathring{\Phi}=\{\mathring{\varphi}^1,\ldots,\mathring{\varphi}^r\}$ and $\mathring{\Psi}=\{\mathring{\psi}^1,\ldots,\mathring{\psi}^s\}$ of $\dLp{2}$ with $s:=r|\det(\dm)|$ such that
\begin{enumerate}
\item[(i)] both $\{\mathring{\Phi}; \Psi\}$ and $\{\mathring{\Phi};\mathring{\Psi}\}$ are $\dm$-framelets in $\dLp{2}$ satisfying \eqref{framelet} (with $\Phi$ and $\Psi$ being replaced by $\mathring{\Phi}$ and $\mathring{\Psi}$, respectively). Moreover, $\mathring{\Psi}$ is a homogeneous $\dm$-framelet in $\dLp{2}$;
\item[(ii)]
$\si(\mathring{\Phi})\subseteq \si(\Phi)$ and $\si(\mathring{\Psi})\subseteq \si_\dm(\Phi)$, i.e., $\wh{\mathring{\varphi}}(\xi)=\wh{\theta}(\xi)\wh{\phi}(\xi)$
and $\wh{\mathring{\psi}}(\dm^\tp\xi)=\wh{\mathring{b}}(\xi)\wh{\phi}(\xi)$ with $\mathring{\varphi}:=(\mathring{\varphi}^1,\ldots,\mathring{\varphi}^r)^\tp$
and $\mathring{\psi}:=(\mathring{\psi}^1,\ldots,\mathring{\psi}^s)^\tp$ for some $r\times r$ matrix $\wh{\theta}$ and some $s\times r$ matrix $\wh{\mathring{b}}$ of $2\pi\dZ$-periodic measurable functions on $\dR$;
\item[(iii)] If $d=1$, $\#\Psi<\infty$ and all the entries of $\wh{b}$ are $2\pi$-periodic trigonometric polynomials, then $\wh{\mathring{b}}$ in item (ii) can be an $s\times r$ matrix of $2\pi$-periodic trigonometric polynomials. If in addition $\phi$ has compact support, then all the elements in $\mathring{\Psi}$ and $\mathring{\psi}$ have compact support.
\end{enumerate}
\end{theorem}

\bp Define $H$ as in \eqref{Psi:H}. By \eqref{refstr}, we have $\si(H)\subset \si(\Phi)$ and $\si(\Psi)\subset \si_{\dm}(\Phi)$. Since $\#\Phi=r$, we have $\len(\si(H))\le \len(\si(\Phi))\le r$ and $\len(\si(\Psi))\le \len(\si_{\dm}(\Phi))\le r|\det(\dm)|$. By Corollary~\ref{cor:hf:1} and Proposition~\ref{prop:hf:1}, there exist $\mathring{\Phi}=\{\mathring{\varphi}^1,\ldots,\mathring{\varphi}^r\}\subset \si(H)\subseteq \si(\Phi)$ and $\mathring{\Psi}=\{\mathring{\psi}^1,\ldots,\mathring{\psi}^s\}\subset
\si(\Psi)\subseteq \si_\dm(\Phi)$ such that item (i) holds. Item (ii) follows directly from $\mathring{\Phi}\subset \si(\Phi)$ and $\mathring{\Psi}\subset \si_{\dm}(\Phi)$.

We now prove item (iii). Note that $\si_{\dm}(\Phi)=\si(\{\phi^\ell_{\dm;k} \setsp \ell=1,\ldots,r, k=0,\ldots,|\dm|-1\})$.
Define $\eta^{(\ell-1)|\dm|+k+1}:=\phi^\ell_{\dm;k}$ for $\ell=1,\ldots,r$ and $k=0,\ldots,|\dm|-1$. Let $\eta:=(\eta^1,\ldots,\eta^s)^\tp$. Then $\si_\dm(\Phi)=\si(\eta)$ and hence the second relation in \eqref{refstr} is equivalent to $\wh{\psi}(\xi)=\Theta(\xi)\wh{\eta}(\xi)$, where $\Theta$ is a unique $\#\Psi\times s$ matrix of $2\pi\dZ$-periodic trigonometric polynomials derived from $\wh{b}$. Using  $\eta^1,\ldots,\eta^s$ instead of
$\varphi^1,\ldots,\varphi^r$ in the proof of Lemma~\ref{lem:hf}, for $f\in \Lp{2}$, by $d=1$, we deduce from \eqref{refstr} that
\be \label{M:identity}
2\pi \sum_{\psi\in \Psi} \sum_{k\in \Z} |\la f, \psi(\cdot-k)\ra|^2=\int_{-\pi}^\pi
[\wh{f},\wh{\eta}](\xi)\ol{\Theta(\xi)}^\tp \Theta(\xi)
[\wh{\eta},\wh{f}](\xi)d\xi.
\ee
Note that $\ol{\Theta(\xi)}^\tp \Theta(\xi)$ is an $s\times s$ positive semidefinite matrix and all its entries are one-dimensional $2\pi$-periodic trigonometric polynomials. By the matrix-valued Fej\'er-Riesz lemma, there exists an $s\times s$ matrix $\wh{v}(\xi)$ of $2\pi$-periodic trigonometric polynomials such that $\ol{\wh{v}(\xi)}^\tp \wh{v}(\xi)=\ol{\Theta(\xi)}^\tp \Theta(\xi)$. Define $\wh{\mathring{\psi}}(\xi):=\wh{v}(\xi)\wh{\eta}(\xi)$.
Then $\si(\mathring{\Psi})\subseteq \si(\eta)=\si_\dm(\Phi)$ and
\eqref{framelet:preserve} holds with $\eta^\ell$ being replaced by $\mathring{\psi}^\ell$ for $\ell=1,\ldots,s$, respectively.
Since $\si(\mathring{\Psi})\subseteq \si(\eta)=\si_\dm(\Phi)$,
we have $\wh{\mathring{\psi}}(\xi)=\wh{\mathring{b}}((\dm^\tp)^{-1}\xi)
\wh{\phi}((\dm^\tp)^{-1}\xi)$ for a unique $s\times r$ matrix of $2\pi$-periodic trigonometric polynomials derived from $\wh{v}$. Consequently, items (i) and (ii) still hold. This proves item (iii).
\ep

If the refinable structure in \eqref{refstr} holds for finite subsets $\Phi$ and $\Psi$ of $\dLp{2}$ and if $\Psi$ is a homogeneous $\dm$-framelet in $\dLp{2}$,
then it has been shown in \cite{ams14,aps16} and \cite{hanbook} for dimension one with $\dm=2$ that there exists a finite subset $\mathring{\Phi}$ with $\#\mathring{\Phi}=\#\Phi$ such that $\{\mathring{\Phi};\Psi\}$ is an $\dm$-framelet in $\dLp{2}$. Therefore,
Item (i) of Theorem~\ref{thm:hf:refstr} generalizes the corresponding results in \cite{ams14,aps16,hanbook} as special cases.

If a homogeneous (tight) $\dm$-framelet $\Psi$ in $\dLp{2}$ is derived from an $r\times 1$ $\dm$-refinable vector function through the refinable structure in \eqref{refstr}, then Theorem~\ref{thm:hf:refstr} tells us that we can always obtain a homogeneous (tight) $\dm$-framelet $\mathring{\Psi}$ and a nonhomogeneous (tight) $\dm$-framelet $\{\mathring{\Phi}; \mathring{\Psi}\}$ in $\dLp{2}$ with $\#\mathring{\Phi}\le r$ and $\#\mathring{\Psi}\le r|\det(\dm)|$.
If in addition the conditions in item (iii) of Theorem~\ref{thm:hf:refstr} are satisfied, then we have a one-dimensional compactly supported homogeneous (tight)
$\dm$-framelet $\mathring{\Psi}$ with $\#\mathring{\Psi}\le r|\dm|$.
Therefore, Theorem~\ref{thm:hf:refstr} is of interest for constructing framelets or tight framelets with a small number of generators.
For example, the compactly supported homogeneous tight $2$-framelet $\Psi$ constructed in Ron and Shen \cite{rs97} from the B-spline of order $m$ has $m$ generators in $\Psi$.
By Theorem~\ref{thm:hf:refstr}, we have a compactly supported homogeneous tight $2$-framelet $\mathring{\Psi}$ with no more than two generators in $\mathring{\Psi}$. As an another example,
the projection method proposed in \cite{han14} (also see \cite{han03}) is a painless method to derive tight $\dm$-framelets $\Psi$ from tensor product tight framelets derived from a scalar refinable function (i.e., $r=1$). The major problem of the projection method in \cite{han14,han03} is that the ratio $\#\Psi/|\det(\dm)|$ is often quite large. Theorem~\ref{thm:hf:refstr} tells us that we can now have a homogeneous tight $\dm$-framelet $\mathring{\Psi}$ so that $\#\mathring{\Psi}=|\det(\dm)|$.

For the special case of homogeneous tight $\dm$-framelets, the following result shows that the condition on the refinable structure in \eqref{refstr} in Theorem~\ref{thm:hf:refstr} can be removed, largely because every nonhomogeneous tight framelet has an intrinsic refinable structure.

\begin{theorem}\label{thm:htf}
Let $\dm$ be a $d\times d$  expansive integer matrix.
Let $\Psi$ be a countable subset of $\dLp{2}$ such that
$\sum_{\psi\in \Psi} \|\psi\|^2_{\dLp{2}}<\infty$.
Define $H$ as in \eqref{Psi:H}
and assume $r:=\len(\si(H))<\infty$. Then the following statements are equivalent:
\begin{enumerate}
\item[(1)] $\Psi$ is a homogeneous tight $\dm$-framelet in $\dLp{2}$ satisfying \eqref{hframelet} with $C_1=C_2=1$.
\item[(2)] There exists a subset $\Phi=\{\varphi^1,\ldots,\varphi^r\}\subset \si(H)$ such that
$\si(\Phi)=\si(H)$ and $\{\Phi; \Psi\}$ is a tight $\dm$-framelet in $\dLp{2}$, i.e.,
\be \label{tightframelet}
\sum_{\varphi\in \Phi}
\sum_{k\in \dZ} |\la f, \varphi(\cdot-k)\ra|^2+
\sum_{j=0}^\infty \sum_{\psi\in \Psi}\sum_{k\in \dZ} |\la f, \psi_{\dm^j;k}\ra|^2=\|f\|_{\dLp{2}}^2,\qquad \forall\, f\in \dLp{2}.
\ee
\item[(3)] There exists a subset $\mathring{\Psi}:=\{\eta^1,\ldots,\eta^s\}\subset \si(\Psi)$ with $s:=r|\det(\dm)|$ such that $\{\Phi;\mathring{\Psi}\}$ is a tight $\dm$-framelet in $\dLp{2}$ and the identity \eqref{framelet:preserve} holds.
\end{enumerate}
\end{theorem}

\bp It follows directly from Proposition~\ref{prop:hf:1} that (1)$\imply$(2).
We now prove (2)$\imply$(3).
Since $\{\Phi;\Psi\}$ is a tight $\dm$-framelet in $\dLp{2}$, by definition, the nonhomogeneous $\dm$-affine system $\AS_0(\Phi;\Psi)$ is a (normalized) tight frame for $\dLp{2}$.
By \cite[Proposition~4]{han12}, every $\AS_J(\Phi;\Psi)$ is a tight frame for $\dLp{2}$, that is,
\be \label{tightframelet:J}
\sum_{\varphi\in \Phi}
\sum_{k\in \dZ} |\la f, \varphi_{\dm^J;k}\ra|^2+
\sum_{j=J}^\infty \sum_{\psi\in \Psi}\sum_{k\in \dZ} |\la f, \psi_{\dm^j;k}\ra|^2=\|f\|_{\dLp{2}}^2,\quad \forall\,f\in \dLp{2}
\ee
for all $J\in \Z$. Considering the difference between $J=0$ and $J=1$, we deduce from \eqref{tightframelet:J} that
\be \label{casecade:tf}
\sum_{\varphi\in \Phi}
\sum_{k\in \dZ} |\la f, \varphi(\cdot-k)\ra|^2
+\sum_{\psi\in \Psi}\sum_{k\in \dZ}
|\la f, \psi(\cdot-k)\ra|^2=
\sum_{\varphi\in \Phi}\sum_{k\in \dZ} |\la f,\varphi_{\dm;k}\ra|^2,\qquad \forall\, f\in \dLp{2}.
\ee
It is straightforward to deduce from the above identity \eqref{casecade:tf} (see \cite{hanrefstr}) that
$\si(\Phi)\cup\si(\Psi)\subseteq \si_\dm(\Phi)$.
By $\len(\si(\Phi))\le r$, we trivially have $\len(\si_{\dm}(\Phi))\le \len(\si(\Phi)) |\det(\dm)|\le r |\det(\dm)|$.
Since $\si(\Psi)\subseteq \si_{\dm}(\Phi)$, we have $\len(\si(\Psi))\le \len(\si_{\dm}(\Phi))\le r|\det(\dm)|$.
Now the conclusion in item (3) follows directly from Corollary~\ref{cor:hf:1}.
This proves (2)$\imply$(3).

Finally, we prove (3)$\imply$(1).
Due to \eqref{framelet:preserve}, item (3) implies that $\{\Phi;\Psi\}$ is a tight $\dm$-framelet in $\dLp{2}$. Since $\dm$ is expansive, now by \cite[Proposition~4]{han12}, we conclude that $\Psi$ must be a homogeneous $\dm$-framelet in $\dLp{2}$. This proves (3)$\imply$(1).
\ep

The main interest of linking a homogeneous tight framelet with a nonhomogeneous tight framelet in Theorem~\ref{thm:htf} lies in that a nonhomogeneous tight framelet always has the refinable structure and is closely related to filter banks.
As proved in \cite{hanrefstr},
all nonhomogeneous tight framelets intrinsically have the refinable structure and are completely characterized in \cite{hanrefstr} (also see \cite[Theorem~4.5.4]{hanbook} for the special case $\dm=2$ in dimension one) through filter banks and the refinable structure as follows.

\begin{theorem} \label{thm:ntf} (\cite{hanrefstr})
Let $\dm$ be a $d\times d$  expansive integer matrix.
Let $\Phi=\{\phi^1,\ldots,\phi^r\}$ and $\Psi=\{\psi^1,\ldots,\psi^s\}$ with $r,s\in \N$ be finite subsets of $\dLp{2}$. Define $\phi:=(\phi^1,\ldots,\phi^r)^\tp$ and $\psi:=(\psi^1,\ldots,\psi^s)^\tp$.
Then $\{\Phi;\Psi\}$ is a tight $\dm$-framelet in $\dLp{2}$ if and only if
\begin{enumerate}
\item[(i)] $\lim_{j\to \infty} \la \|\wh{\phi}((\dm^\tp)^{-j}\cdot)\|_{l_2}^2, h\ra=\la 1, h\ra$ for all compactly supported $\mathscr{C}^\infty$ functions $h$ on $\dR$;
\item[(ii)] There exist an $r\times r$ matrix $\wh{a}$ and an $s\times r$ matrix $\wh{b}$ of $2\pi\dZ$-periodic measurable functions on $\dR$ such that the refinable structure in \eqref{refstr} holds: $\wh{\phi}(\dm^\tp \xi)=\wh{a}(\xi)\wh{\phi}(\xi)$ and $\wh{\psi}(\dm^\tp\xi)=\wh{b}(\xi)\wh{\phi}(\xi)$;
\item[(iii)] $\{\wh{a};\wh{b}\}$ is a generalized tight $\dm$-framelet filter bank, i.e., for all $k\in \dZ$,
\be \label{gtffb}
\begin{split}
&\wh{\phi}(\xi)^\tp \left(\wh{a}(\xi)^\tp \ol{\wh{a}(\xi)}
   +\wh{b}(\xi)^\tp\ol{\wh{b}(\xi)}-I_r\right)\ol{\wh{\phi}(\xi+2\pi k)}=0,\qquad a.e.\, \xi\in \dR,\\
&\wh{\phi}(\xi)^\tp \left(\wh{a}(\xi)^\tp \ol{\wh{a}(\xi+2\pi \omega)}
   +\wh{b}(\xi)^\tp\ol{\wh{b}(\xi+2\pi \omega)}\right)\ol{\wh{\phi}(\xi+2\pi k)}=0,\qquad a.e.\, \xi\in \dR
\end{split}
\ee
for all $\omega\in \dmfc\bs\{0\}$, where $\dmfc:=[(\dm^\tp)^{-1}\dZ]\cap [0,1)^d$.
\end{enumerate}
\end{theorem}

\bp The special case $\dm=2$ has been proved in \cite[Theorem~4.5.4]{hanbook}. The general case is implicitly given in \cite{han12} and can be proved in a similar way as given in \cite[Theorem~4.5.4]{hanbook}.
See \cite{hanrefstr} for a complete detailed proof.
\ep

Theorems~\ref{thm:htf} and~\ref{thm:ntf} show that a homogeneous tight framelet must be intrinsically linked to a nonhomogeneous tight framelet and the refinable structure. Therefore, homogeneous tight framelets can be constructed from filter banks through the refinable structure.

\section{Homogeneous Dual Framelets with the Refinable Structure}

In this section we study the connections of homogeneous dual framelets with nonhomogeneous dual framelets and the refinable structure.

To study homogeneous dual framelets, we need the following result.

\begin{lemma}\label{lem:df:basics}
Let $H$ and $\tilde{H}$ be countable subsets of $\dLp{2}$ such that $\#H=\#\tilde{H}$ with $\sim$ being the bijection between them.
Assume that $\sum_{h\in H} \|h\|_{\dLp{2}}^2<\infty$ and
$\sum_{\tilde{h}\in \tilde{H}} \|\tilde{h}\|_{\dLp{2}}^2<\infty$.
Suppose that there exists a positive constant $C$ such that \eqref{bessel} holds and
\be \label{tilde:bessel}
\sum_{\tilde{h}\in \tilde{H}} \sum_{k\in \dZ} |\la g, \tilde{h}(\cdot-k)\ra|^2\le C\|g\|_{\dLp{2}}^2, \qquad \forall\, g\in \dLp{2}.
\ee
If $r:=\len(\si(H))<\infty$,
then there exist $\eta^1,\ldots,\eta^r\in \si(H)$ and $\tilde{\eta}^1,\ldots,\tilde{\eta}^r \in \si(\tilde{H})$ such that
\begin{enumerate}
\item[(i)] the identity \eqref{basic:identity} holds and
\be \label{df:basic:identity}
\sum_{\ell=1}^r \sum_{k\in \dZ} |\la g, \tilde{\eta}^\ell(\cdot-k)\ra|^2\le \sum_{\tilde{h} \in \tilde{H}} \sum_{k\in \dZ} |\la g, \tilde{h}(\cdot-k)\ra|^2,\qquad \forall\, g\in \dLp{2};
\ee
\item[(ii)] $\si(\{\eta^1,\ldots,\eta^r\})=\si(H)$ and $\si(\{\tilde{\eta}^1,\ldots,\tilde{\eta}^r\})\subseteq \si(\tilde{H})$;
\item[(iii)] The following identity holds:
\be \label{df:identity}
\sum_{\ell=1}^r \sum_{k\in \dZ} \la f, \eta^\ell(\cdot-k)\ra\la \tilde{\eta}^\ell(\cdot-k),g\ra=
\sum_{h \in H} \sum_{k\in \dZ} \la f, h(\cdot-k)\ra \la \tilde{h}(\cdot-k),g\ra,\qquad \forall\, f,g\in \dLp{2}
\ee
with both series converging absolutely.
\end{enumerate}
\end{lemma}

\bp
Let $\varphi^1,\ldots,\varphi^r,
\eta^1,\ldots,\eta^r, u_{h,1},\ldots,u_{h,r}, h\in H$ be as in the proof of Lemma~\ref{lem:hf}.
Then the identity \eqref{basic:identity} holds and $\si(\{\eta^1,\ldots,\eta^r\})=\si(H)$.
Define an $(\#H)\times r$ matrix $B(\xi):=(u_{h,k})_{h\in H, 1\le k\le r}$ and $\varphi:=(\varphi^1,\ldots,\varphi^r)^\tp$. Define $\vec{h}$ to be the column vector function by listing all the elements in $\tilde{H}$ with the same ordering of $H$ as in the matrix $B$.
By \eqref{bessel} and \eqref{tilde:bessel}, for $f,g\in \dLp{2}$, we have
\be \label{fg:identity}
(2\pi)^d \sum_{h\in H} \sum_{k\in \dZ}
\la f, h(\cdot-k)\ra\la \tilde{h}(\cdot-k),g\ra=
\int_{(-\pi,\pi]^d}
[\wh{f}, \wh{\varphi}](\xi) \ol{B(\xi)}^\tp [\wh{\vec{h}},\wh{g}](\xi) d\xi.
\ee
Define $A(\xi):=\ol{B(\xi)}^\tp B(\xi)$, which is the same matrix as in \eqref{matrix:A}. Note that $\eta$ in Lemma~\ref{lem:hf} is defined in \eqref{def:eta:A} as $\wh{\eta}(\xi):=A^{1/2}(\xi) \wh{\varphi}(\xi)$ with
\[
A^{1/2}(\xi)=U(\xi)\mbox{diag}\left(
\sqrt{\gl_1(\xi)},
\ldots,\sqrt{\gl_r(\xi)}\right)\ol{U(\xi)}^\tp,
\]
where $U$ is an $r\times r$ unitary matrix of measurable functions and $\gl_j(\xi)\ge 0$ for all $j=1,\ldots,r$. Let $D(\xi)$ be the pseudoinverse of $A^{1/2}(\xi)$, i.e., $D(\xi):=U(\xi)\mbox{diag}\left(1/\sqrt{\gl_1(\xi)},
\ldots,1/\sqrt{\gl_r(\xi)}\right)\ol{U(\xi)}^\tp$, under the convention that $1/\sqrt{\gl_j(\xi)}:=0$ if $\gl_j(\xi)=0$.
We now claim that
\be \label{ABC}
\ol{B(\xi)}^\tp=A^{1/2}(\xi) D(\xi) \ol{B(\xi)}^\tp.
\ee
As in Theorem~\ref{cor:svd}, it is trivial to observe that $A^{1/2}(\xi)D(\xi)=U(\xi) \mbox{diag}\left(\mbox{sgn}(\gl_1(\xi)),\ldots, \mbox{sgn}(\gl_r(\xi))\right) \ol{U(\xi)}^\tp$.
Therefore, \eqref{ABC} is equivalent to
\be \label{BB}
F(\xi)=
\mbox{diag}\left(\mbox{sgn}(\gl_1(\xi)),\ldots, \mbox{sgn}(\gl_r(\xi))\right)F(\xi) \quad \mbox{with}\quad F(\xi):=\ol{U(\xi)}^\tp \ol{B(\xi)}^\tp.
\ee
On the other hand, by $\ol{B(\xi)}^\tp B(\xi)=A(\xi)=A^{1/2}(\xi)A^{1/2}(\xi)=U(\xi) \mbox{diag}\left(\gl_1(\xi),\ldots,\gl_r(\xi)\right) \ol{U(\xi)}^\tp$, since $U$ is a unitary matrix, we have
\[
F(\xi)\ol{F(\xi)}^\tp=
\ol{U(\xi)}^\tp \ol{B(\xi)}^\tp B(\xi) U(\xi)=
\mbox{diag}\left(\gl_1(\xi),\ldots,\gl_r(\xi)\right),
\]
which implies that if $\gl_j(\xi)=0$, then the $j$th row of $F(\xi)$ must be the zero vector. Consequently, the identity in \eqref{BB} trivially holds and therefore, \eqref{ABC} holds.

By \eqref{ABC}, we have
$\ol{B(\xi)}^\tp=A^{1/2}(\xi)V(\xi)$ with $V(\xi):=D(\xi) \ol{B(\xi)}^\tp$.
Since $D(\xi)$ is the pseudoinverse of $A^{1/2}(\xi)$, we deduce that
\[
V(\xi)\ol{V(\xi)}^\tp
=D(\xi)\ol{B(\xi)}^\tp B(\xi)\ol{D(\xi)}^\tp=
D(\xi) A(\xi) \ol{D(\xi)}^\tp=
\mbox{diag}\left(\mbox{sgn}(\gl_1(\xi)),\ldots, \mbox{sgn}(\gl_r(\xi))\right).
\]
That is, we proved that
\be \label{V:1:0}
\mbox{The nonzero rows of the $r\times (\#H)$ matrix $V(\xi)$ form an orthonormal system in $l_2(\C^{\#H})$}.
\ee
Define $\tilde{\eta}=(\tilde{\eta}^1,\ldots,\tilde{\eta}^r)^\tp$ through $\wh{\tilde{\eta}}(\xi):=V(\xi)\wh{\vec{h}}(\xi)$.
By our assumption $\sum_{\tilde{h}\in \tilde{H}} \|\tilde{h}\|_{\dLp{2}}^2<\infty$, using the Cauchy-Schwarz inequality, we deduce directly from \eqref{V:1:0} that all the entries in $\tilde{\eta}$ belong to $\dLp{2}$. Consequently, it follows from $\wh{\tilde{\eta}}(\xi)=V(\xi)\wh{\vec{h}}(\xi)$ that $\si(\{\tilde{\eta}^1,\ldots,\tilde{\eta}^r\})\subseteq \si(\tilde{H})$.

By $\ol{B(\xi)}^\tp=A^{1/2}(\xi)V(\xi)$ and $\wh{\eta}(\xi)=A^{1/2}(\xi)\wh{\varphi}(\xi)$, we conclude that
\begin{align*}
\int_{(-\pi,\pi]^d}
[\wh{f}, \wh{\varphi}](\xi) \ol{B(\xi)}^\tp  [\wh{\vec{h}},\wh{g}](\xi) d\xi
&=\int_{(-\pi,\pi]^d}
[\wh{f}, \wh{\varphi}](\xi) A^{1/2}(\xi) V(\xi)[\wh{\vec{h}},\wh{g}](\xi) d\xi\\
&=\int_{(-\pi,\pi]^d}
[\wh{f}, A^{1/2}\wh{\varphi}](\xi) [V\wh{\vec{h}},\wh{g}](\xi) d\xi
=\int_{(-\pi,\pi]^d}
[\wh{f}, \wh{\eta}](\xi) [\wh{\tilde{\eta}},\wh{g}](\xi) d\xi,
\end{align*}
from which and \eqref{fg:identity}, we conclude that \eqref{df:identity} holds.

We complete the proof by proving \eqref{df:basic:identity}.
Similarly, by $\wh{\tilde{\eta}}(\xi)=V(\xi) \wh{\vec{h}}(\xi)$, for $g\in \dLp{2}$, we have
\begin{align*}
(2\pi)^d \sum_{\ell=1}^r \sum_{k\in \dZ}
|\la g, \tilde{\eta}^\ell(\cdot-k)\ra|^2
&=\int_{(-\pi,\pi]^d}
\|[\wh{\tilde{\eta}},\wh{g}](\xi)\|_{l_2}^2 d\xi
=\int_{(-\pi,\pi]^d}
\| V(\xi) [\wh{\vec{h}},\wh{g}](\xi)\|_{l_2}^2 d\xi\\
&\le \int_{(-\pi,\pi]^d}
\|[\wh{\vec{h}},\wh{g}](\xi)\|_{l_2}^2 d\xi
=(2\pi)^d \sum_{\tilde{h}\in \tilde{H}}
\sum_{k\in \dZ} |\la g,\tilde{h}(\cdot-k)\ra|^2,
\end{align*}
where we used \eqref{V:1:0} to prove
$\| V(\xi) [\wh{\vec{h}},\wh{g}](\xi)\|_{l_2}
\le \|[\wh{\vec{h}},\wh{g}](\xi)\|_{l_2}$.
This proves the inequality \eqref{df:basic:identity}.
\ep

We remark here that Lemma~\ref{lem:df:basics} can be also proved using singular value decomposition for matrices of measurable functions in Corollary~\ref{cor:svd}.
As a direct consequence of Lemma~\ref{lem:df:basics}, we have

\begin{cor}\label{cor:hdf:1}
Let $\dm$ be a $d\times d$ invertible real-valued matrix.
Let $\Psi$ and $\tilde{\Psi}$ be countable subsets of $\dLp{2}$ with $\sim$ being the bijection between them such that $\sum_{\psi\in \Psi}\|\psi\|_{\dLp{2}}^2<\infty$ and
$\sum_{\tilde{\psi}\in \tilde{\Psi}}\|\tilde{\psi}\|_{\dLp{2}}^2<\infty$.
Suppose that $s:=\min(\len(\si(\Psi)),\len(\si(\tilde{\Psi})))<\infty$.
If $\{\tilde{\Psi},\Psi\}$ is a homogeneous dual $\dm$-framelet in $\dLp{2}$, then there exist subsets $H:=\{\eta^1,\ldots,\eta^s\}\subset \si(\Psi)$ and $\tilde{H}=\{\tilde{\eta}^1,\ldots,\tilde{\eta}^s\}\subset
\si(\tilde{\Psi})$ such that
$(\tilde{H},H)$ is a homogeneous dual $\dm$-framelet in $\dLp{2}$.
\end{cor}

\bp Without loss of generality, we assume $s=\len(\si(\Psi))$. Let $\eta^1,\ldots,\eta^s,\tilde{\eta}^1,\ldots,\tilde{\eta}^s$ be constructed in Lemma~\ref{lem:df:basics}.
Since $\AS(\Psi)$ is a frame in $\dLp{2}$ satisfying \eqref{hframelet},
by Lemma~\ref{lem:df:basics},
$\AS(H)$ also satisfies \eqref{hframelet} with $\Psi$ being replaced by $H$. Therefore, $H$ is a homogeneous $\dm$-framelet in $\dLp{2}$ and $\si(H)=\si(\Psi)$.
By \eqref{df:basic:identity} with $r=s$, similarly we conclude that $\AS(\tilde{H})$ is a Bessel sequence in $\dLp{2}$.
It also follows from the identity \eqref{df:identity} with $r=s$ and \eqref{hdf} that
\[
\sum_{j\in \Z} \sum_{\ell=1}^s \sum_{k\in \dZ} \la f, \eta^\ell_{\dm^j;k}\ra\la \tilde{\eta}^\ell_{\dm^j;k},g\ra=\la f,g\ra,\qquad \forall\, f,g\in \dLp{2}.
\]
Since $\AS(H)$ is a frame of $\dLp{2}$ and $\AS(\tilde{H})$ is a Bessel sequence in $\dLp{2}$, using the Cauchy-Schwarz inequality, we conclude from the above identity that
\[
|\la f, g\ra|^2\le \left(\sum_{j\in \Z} \sum_{\ell=1}^s \sum_{k\in \dZ} |\la f, \eta^\ell_{\dm^j;k}\ra|^2\right)
\left(\sum_{j\in \Z} \sum_{\ell=1}^s \sum_{k\in \dZ} |\la g, \tilde{\eta}^\ell_{\dm^j;k}\ra|^2\right)
\le C_2\|f\|_{\dLp{2}}^2 \left(\sum_{j\in \Z} \sum_{\ell=1}^s \sum_{k\in \dZ} |\la g, \tilde{\eta}^\ell_{\dm^j;k}\ra|^2\right).
\]
Consequently, we deduce from the above inequality that
$\sum_{j\in \Z} \sum_{\ell=1}^s \sum_{k\in \dZ} |\la g, \tilde{\eta}^\ell_{\dm^j;k}\ra|^2
\ge \frac{1}{C_2}\|g\|_{\dLp{2}}^2$.
Thus, $\AS(\tilde{H})$ must be a frame of $\dLp{2}$, that is, $\tilde{H}$ is a homogeneous $\dm$-framelet in $\dLp{2}$.
Therefore, we proved that $(\tilde{H},H)$ is a homogeneous dual $\dm$-framelet in $\dLp{2}$.
\ep

Let $\Phi,\Psi,\tilde{\Phi}, \tilde{\Psi}$ be subsets of $\dLp{2}$.
We say that $(\{\tilde{\Phi}; \tilde{\Psi}\}, \{\Phi;\Psi\})$ is \emph{a dual $\dm$-framelet} in $\dLp{2}$ if
$(\AS_0(\tilde{\Phi}; \tilde{\Psi}), \AS_0(\Phi;\Psi))$ is a pair of dual frames in $\dLp{2}$, that is,
(i) each of $\AS_0(\Phi;\Psi)$ and $\AS_0(\tilde{\Phi}; \tilde{\Psi})$ is a frame for $\dLp{2}$, and (ii) the following identity holds
\be \label{df}
\sum_{\varphi\in \Phi}\sum_{k\in \dZ}
\la f, \varphi(\cdot-k)\ra\la \tilde{\varphi}(\cdot-k), g\ra+
\sum_{j=0}^\infty \sum_{\psi\in \Psi} \sum_{k\in \dZ} \la f, \psi_{\dm^j;k}\ra\la \tilde{\psi}_{\dm^j;k},g\ra=\la f,g\ra,\qquad \forall\, f,g\in \dLp{2}.
\ee
We now connect a homogeneous dual framelet with a nonhomogeneous dual framelet as follows.

\begin{prop}\label{prop:hdf}
Let $\dm$ be a $d\times d$ expansive integer matrix.
Let $\Psi$ and $\tilde{\Psi}$ be countable subsets of $\dLp{2}$ with $\sim$ being the bijection between them such that $\sum_{\psi\in \Psi}\|\psi\|_{\dLp{2}}^2<\infty$ and
$\sum_{\tilde{\psi}\in \tilde{\Psi}}\|\tilde{\psi}\|_{\dLp{2}}^2<\infty$.
Define two subsets $H,\tilde{H}$ as follows:
\be \label{Psi:H:tH}
H:=\{ |\det(\dm)|^{-j} \psi(\dm^{-j}\cdot) \setsp j\in \N, \psi\in \Psi\},\quad \tilde{H}:=\{ |\det(\dm)|^{-j} \tilde{\psi}(\dm^{-j}\cdot) \setsp j\in \N, \tilde{\psi}\in \tilde{\Psi}\}.
\ee
Suppose $r:=\len(\si(H))<\infty$.
Then $(\tilde{\Psi},\Psi)$ is a homogeneous dual $\dm$-framelet in $\dLp{2}$ satisfying \eqref{hdf} if and only if there exist $\Phi:=\{\varphi^1,\ldots,\varphi^r\}\subset \si(H)$ and $\tilde{\Phi}:=\{\tilde{\varphi}^1,\ldots,
\tilde{\varphi}^r\}\subset \si(\tilde{H})$ such that
$\si(\Phi)=\si(H)$, \eqref{H:Phi} holds, and
$(\{\tilde{\Phi}; \tilde{\Psi}\}, \{\Phi;\Psi\})$ is a dual $\dm$-framelet in $\dLp{2}$ satisfying \eqref{df}.
\end{prop}

\bp Sufficiency ($\Rightarrow$). By \cite[Theorem~5.5]{rs97} and \cite[Theorem~2]{css98} for a finite set $\Psi$ and \cite[Theorem~4.3.4]{hanbook} for a countable set $\Psi$,
$(\{\tilde{H};\tilde{\Psi}\}, \{H;\Psi\})$ is a dual $\dm$-framelet in $\dLp{2}$.
Now the claim follows directly from Lemma~\ref{lem:df:basics}.

Necessity ($\Leftarrow$).
If $(\{\tilde{\Phi};\tilde{\Psi}\}, \{\Phi;\Psi\})$ is a dual $\dm$-framelet in $\dLp{2}$, since $\dm$ is expansive, then it has been proved in \cite[Proposition~5]{han12} that $(\tilde{\Psi},\Psi)$ must be a homogeneous dual $\dm$-framelet in $\dLp{2}$.
This completes the proof.
\ep

For homogeneous dual framelets with the refinable structure in \eqref{refstr}, we have

\begin{theorem}\label{thm:hdf:refstr}
Let $\dm$ be a $d\times d$ expansive integer matrix.
Let $\Phi=\{\phi^1,\ldots,\phi^r\}$ be a finite subset of $\dLp{2}$ and $\Psi$ be a countable subset of $\dLp{2}$ such that $\sum_{\psi\in \Psi} \|\psi\|_{\dLp{2}}^2<\infty$ and $\si(\Phi)\cup\si(\Psi)\subseteq \si_\dm(\Phi)$ (i.e.,
the refinable structure in \eqref{refstr} holds with $\phi:=(\phi^1,\ldots,\phi^r)^\tp$
and $\psi$ being a column vector by listing all elements in $\Psi$ for some matrices $\wh{a}$ and $\wh{b}$ of $2\pi\dZ$-periodic measurable functions).
Let $\tilde{\Psi}$ be a countable subset of $\dLp{2}$ such that $\sum_{\tilde{\psi}\in \tilde{\Psi}} \|\tilde{\psi}\|_{\dLp{2}}^2<\infty$ and $\#\tilde{\Psi}=\#\Psi$ with $\sim$ being the bijection between them.
If $(\tilde{\Psi},\Psi)$ is a homogeneous dual $\dm$-framelet in $\dLp{2}$,
then there exist subsets
$\mathring{\Phi}=\{\mathring{\varphi}^1,\ldots,\mathring{\varphi}^r\}\subset \si(\Phi)$, $\mathring{\Psi}=\{\mathring{\psi}^1,\ldots,\mathring{\psi}^s\}\subset \si_{\dm}(\Phi)$
with $s:=r|\det(\dm)|$, and subsets
\[
\tilde{\mathring{\Phi}}=\{\tilde{\mathring{\varphi}}^1,
\ldots,\tilde{\mathring{\varphi}}^r\}\subset \si(\tilde{H}),\quad \tilde{\mathring{\Psi}}=\{\tilde{\mathring{\psi}}^1,\ldots,\tilde{\mathring{\psi}}^s\}\subset \si(\tilde{\Psi})
\]
with $\tilde{H}$ being defined in \eqref{Psi:H:tH} such that
\begin{enumerate}
\item[(i)] both $(\{\tilde{\mathring{\Phi}}; \tilde{\Psi}\},\{\mathring{\Phi}; \Psi\})$ and $(\{\tilde{\mathring{\Phi}};\tilde{\mathring{\Psi}}\},
    \{\mathring{\Phi};\mathring{\Psi}\})$ are dual $\dm$-framelets in $\dLp{2}$. Moreover, $(\tilde{\mathring{\Psi}},\mathring{\Psi})$ is a homogeneous dual $\dm$-framelet in $\dLp{2}$ with $\#\tilde{\mathring{\Psi}}=\#\mathring{\Psi}\le r|\det(\dm)|$;
\item[(ii)] $\si(\mathring{\Phi})\subseteq \si(\Phi)$ and $\si(\mathring{\Psi})\subseteq \si_\dm(\Phi)$, i.e., $\wh{\mathring{\varphi}}(\xi)=\wh{\theta}(\xi)\wh{\phi}(\xi)$
and $\wh{\mathring{\psi}}(\dm^\tp\xi)=\wh{\mathring{b}}(\xi)\wh{\phi}(\xi)$ with $\mathring{\varphi}:=(\mathring{\varphi}^1,\ldots,\mathring{\varphi}^r)^\tp$
and $\mathring{\psi}:=(\mathring{\psi}^1,\ldots,\mathring{\psi}^s)^\tp$ for some $r\times r$ matrix $\wh{\theta}$ and $s\times r$ matrix $\wh{\mathring{b}}$ of $2\pi\dZ$-periodic measurable functions on $\dR$.
\end{enumerate}
If in addition
$\si(\tilde{\Phi})\cup\si(\tilde{\Psi})\subseteq \si_\dm(\tilde{\Phi})$ for some finite subset $\tilde{\Phi}=\{\tilde{\phi}^1,\ldots,\tilde{\phi}^{\tilde{r}}\}$ of $\dLp{2}$, then
%
%
\begin{enumerate}
\item[(iii)] $\si(\tilde{\mathring{\Phi}})\subseteq \si(\tilde{\Phi})$ and
$\si(\tilde{\mathring{\Psi}})\subseteq \si_\dm(\tilde{\Phi})$, i.e.,
$\wh{\tilde{\mathring{\varphi}}}(\xi)
=\wh{\tilde{\theta}}(\xi)\wh{\tilde{\phi}}(\xi)$
and $\wh{\tilde{\mathring{\psi}}}(\dm^\tp\xi)=
\wh{\tilde{\mathring{b}}}(\xi)\wh{\tilde{\phi}}(\xi)$ for some $\tilde{r}\times \tilde{r}$ matrix $\wh{\tilde{\theta}}$ and $s\times \tilde{r}$ matrix $\wh{\tilde{\mathring{b}}}$ of $2\pi\dZ$-periodic measurable functions on $\dR$, where $\tilde{\phi}:=(\tilde{\phi}^1,\ldots,\tilde{\phi}^{\tilde{r}})^\tp$ and $\tilde{\mathring{\psi}}:=(\tilde{\mathring{\psi}}^1,\ldots,\tilde{\mathring{\psi}}^{\tilde{r}})^\tp$.
\end{enumerate}
\end{theorem}

\bp Let $H$ and $\tilde{H}$ be defined in \eqref{Psi:H:tH}. By our assumption $\si(\Phi)\cup \si(\Psi)\subseteq \si_\dm(\Phi)$, we conclude that $\si(H)\subseteq \si(\Phi)$ and $\si(\Psi)\subseteq \si_\dm(\Phi)$.
By the same proof of Theorem~\ref{thm:hf:refstr},
the claim is a direct consequence of Corollary~\ref{cor:hdf:1} and Proposition~\ref{prop:hdf}.
\ep

Under the condition that both $\Psi$ and $\tilde{\Psi}$ are obtained from refinable functions through the refinable structure, item (i) of
Theorem~\ref{thm:hdf:refstr} is known in 
\cite{ams14,aps16,hanbook} for the existence of a dual $\dm$-framelet $(\{\tilde{\mathring{\Phi}}; \tilde{\Psi}\},\{\mathring{\Phi}; \Psi\})$. Therefore, Theorem~\ref{thm:hdf:refstr} generalizes \cite{ams14,aps16,hanbook} under a weaker assumption.

The interest of linking a homogeneous dual framelet with the refinable structure to a nonhomogeneous dual framelet in Theorem~\ref{thm:hdf:refstr} lies in that nonhomogeneous dual framelets with the refinable structure is well studied and understood in \cite{hanbook,han12,han10} and is closely linked to the oblique extension principle in \cite{dhrs03} (also cf. \cite{chs02}).
Summarizing \cite[Theorems~2, 7, and~16]{han10},
\cite[Theorems~9 and~17]{han12}, and
\cite[Theorems~4.1.10 and~4.3.7]{hanbook}, we explicitly state the following result:

\begin{theorem}
Let $\dm$ be a $d\times d$ expansive integer matrix.
Let $H=\{\eta^1,\ldots,\eta^r\},
\Psi=\{\psi^1,\ldots,\psi^s\}$
and
$\tilde{H}=\{\tilde{\eta}^1,\ldots,\tilde{\eta}^r\},
\tilde{\Psi}=\{\tilde{\psi}^1,\ldots,\tilde{\psi}^s\}$ be finite subsets of $\dLp{2}$.
Define
\[
\eta:=(\eta^1,\ldots,\eta^r)^\tp,
\quad
\psi:=(\psi^1,\ldots,\psi^s)^\tp,\quad
\tilde{\eta}:=(\tilde{\eta}^1,\ldots,\tilde{\eta}^r)^\tp,
\quad
\tilde{\psi}:=(\tilde{\psi}^1,\ldots,\tilde{\psi}^s)^\tp.
\]
Suppose that $\phi, \tilde{\phi}\in (\dLp{2})^r$ are $\dm$-refinable vector functions satisfying
\be \label{reffunc}
\wh{\phi}(\dm^\tp\xi)=\wh{a}(\xi)\wh{\phi}(\xi),\qquad \wh{\tilde{\phi}}(\dm^\tp \xi)=\wh{\tilde{a}}(\xi)
\wh{\tilde{\phi}}(\xi),\qquad a.e.\,\xi\in \dR,
\ee
for some $r\times r$ matrices $\wh{a}$ and $\wh{\tilde{a}}$ of $2\pi\dZ$-periodic measurable functions on $\dR$.
Suppose that all $\eta,\tilde{\eta},\psi, \tilde{\psi}$ are derived from $\phi$ and $\tilde{\phi}$ by $\wh{\eta}(\xi)=\wh{\theta}(\xi)\wh{\phi}(\xi)$,
$\wh{\tilde{\eta}}(\xi)=\wh{\tilde{\theta}}(\xi)\wh{\tilde{\phi}}(\xi)$
and
\be \label{filterbank}
\wh{\psi}(\dm^\tp\xi)=\wh{b}(\xi)\wh{\phi}(\xi),\qquad
\wh{\tilde{\psi}}(\dm^\tp\xi)=\wh{\tilde{b}}(\xi)\wh{\phi}(\xi)
\ee
for some $r\times r$ matrices $\wh{\theta},\wh{\tilde{\theta}}$ and $s\times r$ matrices $\wh{b}, \wh{\tilde{b}}$ of $2\pi\dZ$-periodic measurable functions on $\dR$.
Then $(\{\tilde{H};\tilde{\Psi}\},\{H;\Psi\})$ is a dual $\dm$-framelet in $\dLp{2}$ if and only if
\begin{enumerate}
\item[(i)] $\lim_{j\to \infty}
\la \wh{\tilde{\phi}}((\dm^\tp)^{-j}\cdot)^\tp \wh{\Theta}((\dm^\tp)^{-j}\cdot) \ol{\wh{\phi}((\dm^\tp)^{-j}\cdot)}, h\ra=\la 1,h\ra$ for all compactly supported $\mathscr{C}^\infty$ functions $h$ on $\dR$, where $\wh{\Theta}(\xi):=\wh{\tilde{\theta}}(\xi)^\tp
\ol{\wh{\theta}(\xi)}$;

\item[(ii)] $(\{\wh{\tilde{a}};\wh{\tilde{b}}\}, \{\wh{a};\wh{b}\})_{\wh{\Theta}}$ is a generalized dual $\dm$-framelet filter bank, i.e., for all $k\in \dZ$,
\begin{align}
&\wh{\tilde{\phi}}(\xi)^\tp\Big[ \wh{\tilde{a}}(\xi)^\tp \wh{\Theta}(\dm^\tp \xi)\ol{\wh{a}(\xi)}+\wh{\tilde{b}}(\xi)^\tp \ol{\wh{b}(\xi)}-\wh{\Theta}(\xi)\Big]\ol{\wh{\phi}(\xi+2\pi k)}=0,\qquad a.e.\,\xi\in \dR,\\
&\wh{\tilde{\phi}}(\xi)^\tp \Big[ \wh{\tilde{a}}(\xi)^\tp \wh{\Theta}(\dm^\tp \xi)\ol{\wh{a}(\xi+2\pi \omega)}+\wh{\tilde{b}}(\xi)^\tp \ol{\wh{b}(\xi+2\pi \omega)}\, \Big]\ol{\wh{\phi}(\xi+2\pi k)}=0,\qquad a.e.\,\xi\in \dR
\end{align}
for all $\omega\in \dmfc\bs\{0\}$, where $\dmfc:=[(\dm^\tp)^{-1}\dZ]\cap [0,1)^d$;
\item[(iii)] Both $\AS_0(H;\Psi)$ and $\AS_0(\tilde{H};\tilde{\Psi})$ are Bessel sequences in $\dLp{2}$.
\end{enumerate}
\end{theorem}

\bp Using the same argument as in \cite[Theorems~4.1.9 and 4.1.10]{hanbook} or \cite[Theorem~2]{han10}, by \cite[Corollary~16]{han12}, one can verify that items (i) and (ii) are equivalent to the fact that $(\{\wh{\tilde{H}}; \wh{\tilde{\Psi}}\},\{\wh{H}; \wh{\Psi}\})$ is a pair of frequency-based dual $\dm$-framelet in the distribution space, see \cite[(3.9)]{han12} for its definition. Now the claim follows directly from \cite[Theorem~9]{han12}.
\ep

\section{Link Homogeneous Wavelets to Nonhomogeneous Wavelets with the Refinable Structure}

In this section we shall link homogeneous
wavelets with nonhomogeneous wavelets with the refinable structure.

Note that a homogeneous Riesz $\dm$-wavelet $\Psi$ in $\dLp{2}$ is always a homogeneous $\dm$-framelet in $\dLp{2}$. Therefore, by Proposition~\ref{prop:hf:1} in Section~3, we can link a homogeneous Riesz $\dm$-wavelet $\Psi$ in $\dLp{2}$ to an $\dm$-framelet $\{\Phi;\Psi\}$ in $\dLp{2}$, which however may not be a Riesz $\dm$-wavelet in $\dLp{2}$ any more. As we shall show by an example at the end of this section, this situation indeed can happen.
In comparison with framelets,
because a wavelet does not allow redundancy,  it is more difficult to link a homogeneous wavelet to a nonhomogeneous wavelet if no extra conditions are imposed.
As the major tool for constructing homogeneous wavelets, multiresolution analysis (MRA) is introduced in \cite{mal89,meybook90} and has been extended to the generalized multiresolution analysis (GMRA), e.g., see \cite{bjmp05,bfmp09} and many references therein.
The homogeneous wavelets, which are associated with or obtained from MRA or GMRA, have been extensively investigated in several papers, for example, see \cite{aus95,bow03,han95,hanmsc94,hwbook96,kkl01,lem92,zal99}.
However, to our best knowledge, except the one-dimensional special case in \cite[Subsection~4.5.5]{hanbook}, there are no other results in the literature to explicitly link a homogeneous wavelet to a nonhomogeneous wavelet. Note that a sequence of nonhomogeneous affine systems $\{\AS_J(\Phi;\Psi)\}_{J\in \Z}$ is not always associated with a generalized multiresolution analysis.

Let $\Psi$ and $\tilde{\Psi}$ be subsets of $\dLp{2}$ such that $\#\Psi=\#\tilde{\Psi}$ with $\sim$ being the bijection between them. We say that
$(\tilde{\Psi},\Psi)$ is \emph{a homogeneous biorthogonal $\dm$-wavelet} in $\dLp{2}$ if (i) each of $\Psi$ and $\tilde{\Psi}$ is a homogeneous Riesz $\dm$-wavelet in $\dLp{2}$, that is, each of $\AS(\Psi)$ and $\AS(\tilde{\Psi})$ is a Riesz basis of $\dLp{2}$,
and (ii) $\AS(\Psi)$ and $\AS(\tilde{\Psi})$ are biorthogonal to each other, i.e.,
\be \label{hbw}
\la h, \tilde{h}\ra=1\quad \mbox{and}\quad
\la h,g\ra=0,\quad \mbox{if}\; g \in \AS(\tilde{\Psi})\bs\{\tilde{h}\},
\qquad \forall\; h\in \AS(\Psi),
\ee
where $\sim$ is the induced bijection between $\AS(\Psi)$ and $\AS(\tilde{\Psi})$.
Note that a homogeneous biorthogonal $\dm$-wavelet in $\dLp{2}$ is always a homogeneous dual $\dm$-framelet in $\dLp{2}$. Therefore, Proposition~\ref{prop:hdf} in Section~4 can be applied to obtain a nonhomogeneous dual $\dm$-framelet $(\{\tilde{\Phi};\tilde{\Psi}\},\{\Phi;\Psi\})$, which however may cease to be a biorthogonal $\dm$-wavelet in $\dLp{2}$.

Before proceeding further, we need the following simple facts on Riesz bases of $\dLp{2}$.

\begin{lemma}\label{lem:riesz}
Let $S_1,S_2,T_1,T_2$ be countable subsets of $\dLp{2}$ such that $S_1\cup S_2$ is a Riesz basis for $\dLp{2}$. Define $U_1:=\ol{\mbox{span} (S_1)}^{\|\cdot\|_{\dLp{2}}}$ and $U_2:=\ol{\mbox{span} (S_2)}^{\|\cdot\|_{\dLp{2}}}$.
\begin{enumerate}
\item[(i)] If $T_1$ and $T_2$ are  Riesz bases of $U_1$ and $U_2$, respectively, then $T_1\cup T_2$ is a Riesz basis of $\dLp{2}$.
\item[(ii)] If $T_1\cup T_2$ is a Riesz basis of $\dLp{2}$ such that  $V_1:=\ol{\mbox{span} (T_1)}^{\|\cdot\|_{\dLp{2}}}\subseteq U_1$ and $V_2:=\ol{\mbox{span} (T_2)}^{\|\cdot\|_{\dLp{2}}}\subseteq U_2$, then $V_1=U_1$ and $V_2=U_2$.
\item[(iii)] If $S_1\cup T_2$ is a Riesz basis of $\dLp{2}$, then $P: V_2\rightarrow U_1^\perp$ is bijective,
    $P(T_2)$ must be a Riesz basis of $U_1^\perp$, and $S_1\cup P(T_2)$ must be a Riesz basis of $\dLp{2}$, where $P: \dLp{2}\rightarrow U_1^\perp$ is the orthogonal projection.
\end{enumerate}
\end{lemma}

\bp
Note that every subset of a Riesz sequence is also a Riesz sequence with the same (or even better) lower and upper Riesz bounds.
Therefore, $S_1$ is a Riesz basis for $U_1$ and $S_2$ is a Riesz basis for $U_2$. Since $S_1\cup S_2$ is a Riesz basis of $\dLp{2}$ with a lower Riesz bound $C_3$ and an upper Riesz bound $C_4$, we have $U_1+U_2=\dLp{2}$ and every $f\in \dLp{2}$ can be uniquely decomposed as $f=f_1+f_2$ with $f_1\in U_1$ and $f_2\in U_2$ such that
\[
\frac{C_3}{C_4}(\|f_1\|_{\dLp{2}}^2+\|f_2\|_{\dLp{2}}^2)
\le \|f_1+f_2\|_{\dLp{2}}^2\le \frac{C_4}{C_3}(\|f_1\|_{\dLp{2}}^2+\|f_2\|_{\dLp{2}}^2),\qquad
\forall\; f_1\in U_1, f_2\in U_2.
\]
To prove item (i), it follows trivially from the above inequalities that $T_1\cup T_2$ must a Riesz sequence. Since $U_1+U_2=\dLp{2}$, we conclude that $T_1\cup T_2$ is a Riesz basis for $\dLp{2}$.

To prove item (ii), since $T_1\cup T_2$ is a Riesz basis for $\dLp{2}$, we must have $V_1+V_2=\dLp{2}$. If there exists $f\in U_1\bs V_1$, then we can uniquely write $f=f_1+f_2$ with $f_1\in V_1 \subseteq U_1$ and $f_2\in V_2\subseteq U_2$.
However, $f=f+0$ is another decomposition of $f$ with $f\in U_1$ and $0\in U_2$. By the uniqueness of the decomposition, we must have $f=f_1\in V_1$, which is a contradiction to $f\in U_1\bs V_1$. Therefore, $V_1=U_1$. Similarly, we can prove $V_2=U_2$.

To prove item (iii), we first show that $P: V_2 \rightarrow U_1^\perp$ is a bijection. If $Pf=0$ for some $f\in V_2$, since $P$ is the orthogonal projection to $U_1^\perp$, then $f\in U_1$. Since $S_1\cup T_2$ is a Riesz basis for $\dLp{2}$, we must have $f\in U_1\cap V_2=\{0\}$ by $f\in V_2$ and $f\in U_1$. Hence, $f=0$ and $P$ is injective.
On the other hand, for $f\in U_1^\perp$, since $S_1\cup T_2$ is a Riesz basis for $\dLp{2}$, we have $f=g+h$ with $g\in U_1$ and $h\in V_2$. Therefore, we have
$f=g+h=g+(h-Ph)+Ph$ and hence, $f-Ph=g+(h-Ph)$. Since $f\in U_1^\perp$ and $Ph\in U_1^\perp$, we have $f-Ph\in U_1^\perp$. However, $h-Ph\in U_1$ and $g\in U_1$. Therefore, $f-Ph=g+(h-Ph)\in U_1$. This shows that $f-Ph\in U_1^\perp\cap U_1=\{0\}$. So, $f=Ph\in PV_2$. This proves that $P$ is surjective. Hence, $P: V_2\rightarrow U_1^\perp$ is a bijection. As an orthogonal projection, the operator $P$ has operator norm one. Because both $V_2$ and $U_1^\perp$ are closed linear spaces of $\dLp{2}$, by the Open Mapping Theorem, we conclude that $P^{-1}:U_1^\perp \rightarrow V_2$ is a bounded operator. Since $T_2$ must be a Riesz basis of $V_2$ and $P$ is an isomorphism between $V_2$ and $U_1^\perp$, $P(T_2)$ must be a Riesz basis of $U_1^\perp$.
Since $S_1$ is a Riesz basis of $U_1$, $S_1\cup P(T_2)$ must be a Riesz basis of $U_1\oplus U_1^\perp=\dLp{2}$.
\ep

We also need the following known simple fact (see \cite{kkl01,ltw01}) for which we provide a proof here.

\begin{lemma}\label{lem:V}
Let $\dm$ be a $d\times d$ invertible integer matrix.
Suppose that $\Psi\subset \dLp{2}$ is a homogeneous Riesz $\dm$-wavelet in $\dLp{2}$. Then there exists $\tilde{\Psi}\subset \dLp{2}$ such that $(\tilde{\Psi},\Psi)$ is a homogeneous biorthogonal $\dm$-wavelet in $\dLp{2}$ if and only if
\be \label{Vminus:cond}
\V_-(\Psi):=\ol{\mbox{span}\AS_-(\Psi)}^{\|\cdot\|_{\dLp{2}}}
\quad \mbox{with}\quad
\AS_-(\Psi):=\{\psi_{\dm^{-j};k} \setsp k\in \dZ, j\in \N, \psi\in \Psi\}
\ee
is shift-invariant. 
%
%
\end{lemma}

\bp
The claim is known in \cite[Theorem~2.8]{kkl01} for dimension $d=1$, $\dm=2$ and $\#\Psi=1$, and
in \cite{ltw01} for high dimensions.
For completeness, we provide a proof here.

Necessity ($\Rightarrow$). Define
\be \label{Vplus:cond}
\V_+(\Psi):=\ol{\mbox{span}\AS_+(\Psi)}^{\|\cdot\|_{\dLp{2}}}
\quad \mbox{with}\quad
\AS_+(\Psi):=\{\psi_{\dm^{j};k} \setsp k\in \dZ, j\in \N\cup\{0\}, \psi\in \Psi\}.
\ee
If $f\in \V_+(\tilde{\Psi})^\perp$, then
$\la f,g\ra=0$ for all $g\in \AS_+(\tilde{\Psi})$ and consequently,
the representation $f=\sum_{j\in \Z} \sum_{\psi\in \Psi}\sum_{k\in \dZ}\la f, \tilde{\psi}_{\dm^j;k}\ra \psi_{\dm^j;k}$
directly tells us that
$f\in \V_-(\Psi)$. Hence, $\V_+(\tilde{\Psi})^\perp\subseteq
\V_-(\Psi)$.
Conversely, the biorthogonality relation between $\AS(\Psi)$ and $\AS(\tilde{\Psi})$ trivially tells us that
$\AS_-(\Psi)\subset \V_+(\tilde{\Psi})^\perp$.
Thus, $\V_-(\Psi)=\ol{\mbox{span}(\AS_-(\Psi))}^{\|\cdot\|_{\dLp{2}}}
\subseteq \V_+(\tilde{\Psi})^\perp$.
In conclusion, we proved $\V_-(\Psi)=\V_+(\tilde{\Psi})^\perp$.
Since $V_+(\tilde{\Psi})$ is obviously shift-invariant, the space $V_-(\Psi)=V_+(\tilde{\Psi})^\perp$ must be shift-invariant.

Sufficiency ($\Leftarrow$).
Suppose that $\V_-(\Psi)$ is shift-invariant.
For each fixed $\psi\in \Psi$,
since $\AS(\Psi)$ is a Riesz basis of $\dLp{2}$, there is a unique element
$\tilde{\psi}\in \dLp{2}$ satisfying
\be \label{dual:riesz}
\la \tilde{\psi},\psi\ra=1 \quad \mbox{and}\quad \la \tilde{\psi},h\ra=0,\qquad \forall\, h\in \AS(\Psi)\bs \{\psi\}.
\ee
Therefore, we constructed a set $\tilde{\Psi}$ such that \eqref{dual:riesz} is satisfied for every $\psi\in \Psi$.
Let $\tilde{\psi}\in \tilde{\Psi}$ and $\eta\in \Psi$. Let $j,j'\in \Z$ and $k,k'\in \dZ$. If $j\le j'$, by \eqref{dual:riesz}, then we have $\dm^{j'-j}k\in \dZ$ and
\[
\la \tilde{\psi}_{\dm^j;k},\eta_{\dm^{j'};k'}\ra=
\la \tilde{\psi},\eta_{\dm^{j'-j};k'-\dm^{j'-j}k}\ra=
\begin{cases}
1, &\text{if $\tilde{\eta}=\tilde{\psi}, j=j',k=k'$},\\
0, &\text{otherwise}.\end{cases}
\]
If $j>j'$, then $\la \tilde{\psi}_{\dm^j;k},\eta_{\dm^{j'};k'}\ra
=\la \tilde{\psi},\eta_{\dm^{j'-j};k'}(\cdot+k)\ra=0$, since $\V_-(\Psi)$ is shift-invariant.
This proves that $\AS(\tilde{\Psi})$ is biorthogonal to $\AS(\Psi)$.
Hence, $\AS(\tilde{\Psi})$ must be the unique dual Riesz basis of
$\AS(\Psi)$.
Consequently, $(\tilde{\Psi},\Psi)$ is a homogeneous biorthogonal $\dm$-wavelet in $\dLp{2}$.
\ep


Recall that $(\{\tilde{\Phi};\tilde{\Psi}\},\{\Phi;\Psi\})$ is \emph{a biorthogonal $\dm$-wavelet} in $\dLp{2}$ if (i) each of $\{\Phi;\Psi\}$ and $\{\tilde{\Phi};\tilde{\Psi}\}$ is a Riesz $\dm$-wavelet in $\dLp{2}$, and (ii) $\AS_0(\Phi;\Psi)$ and $\AS_0(\tilde{\Phi};\tilde{\Psi})$ are biorthogonal to each other.
If $\dm$ is a $d\times d$ invertible integer matrix and $(\Phi;\Psi)$ is a Riesz $\dm$-wavelet in $\dLp{2}$, then it is known in \cite[Theorem~8]{han12} that there exist $\tilde{\Phi},\tilde{\Psi}\subseteq \dLp{2}$ such that $(\{\tilde{\Phi};\tilde{\Psi}\},\{\Phi;\Psi\})$ is a biorthogonal $\dm$-wavelet in $\dLp{2}$ if and only if $\Phi$ and $\Psi$ must have the refinable structure $\si(\Phi)\cup\si(\Psi)\subseteq \si_\dm(\Phi)$, that is,
\be \label{refstr:wavelet}
\wh{\phi}(\dm^\tp\xi)=\wh{a}(\xi)\wh{\phi}(\xi),\quad
\wh{\psi}(\dm^\tp\xi)=\wh{b}(\xi)\wh{\phi}(\xi),\qquad a.e.\,\xi\in \dR,
\ee
where $\phi$ and $\psi$ are column vector functions by listing the elements in $\Phi$ and $\Psi$, respectively, and $\wh{a}$ and $\wh{b}$ are matrices of $2\pi\dZ$-periodic measurable functions on $\dR$.
It is also known in \cite[Theorems~7 and~8]{han12} that if $(\{\tilde{\Phi};\tilde{\Psi}\},\{\Phi;\Psi\})$ is a biorthogonal $\dm$-wavelet in $\dLp{2}$, then $(\AS_J(\tilde{\Phi};\tilde{\Psi}), \AS_J(\Phi;\Psi))$ is a pair of biorthogonal Riesz bases in $\dLp{2}$ for all $J\in \Z$. If in addition $\dm$ is an expansive integer matrix and $\sum_{\phi\in \Phi}(\|\phi\|_{\dLp{2}}^2
+\|\tilde{\phi}\|_{\dLp{2}}^2)<\infty$, then
$(\tilde{\Psi},\Psi)$ must be a homogeneous biorthogonal $\dm$-wavelet in $\dLp{2}$; More importantly, the biorthogonal $\dm$-wavelet $(\{\tilde{\Phi};\tilde{\Psi}\},\{\Phi;\Psi\})$ must have the refinable structure, see Theorem~\ref{thm:nbw} below for details.

We now study a homogeneous Riesz $\dm$-wavelet $\Psi$ in $\dLp{2}$ under the condition that $\V_-(\Psi)$ in \eqref{Vminus:cond} is shift-invariant.

\begin{theorem}\label{thm:hw}
Let $\dm$ be a $d\times d$ expansive integer matrix and $\Psi=\{\psi^1,\ldots,\psi^s\}$ be a finite subset of $\dLp{2}$. Then the following statements are equivalent:
\begin{enumerate}
\item[(i)] $\Psi$ is a homogeneous Riesz $\dm$-wavelet in $\dLp{2}$, $\V_-(\Psi)$ is shift-invariant, and $\dim_{\V_-(\Psi)}(\xi)=r_1$ for almost every $\xi\in \dR$ for some constant $r_1$.

\item[(ii)] There exists a subset $\tilde{\Psi}=\{\tilde{\psi}^1,\ldots,\tilde{\psi}^s\}\subset \dLp{2}$, which is uniquely determined by the identities in \eqref{dual:riesz},
     such that $(\tilde{\Psi},\Psi)$ is a homogeneous biorthogonal $\dm$-wavelet in $\dLp{2}$, and
\be \label{dim:Psi}
\sum_{j=1}^\infty\sum_{\ell=1}^s [\wh{\psi^\ell}((\dm^\tp)^{j}\cdot), \wh{\tilde{\psi}^\ell}((\dm^\tp)^j\cdot)](\xi)=r_2,\quad a.e.\, \xi\in \dR
\ee
for some constant $r_2$.

\item[(iii)] There exists a subset $\Phi=\{\varphi^1,\ldots,\varphi^r\} \subset \dLp{2}$ such that $\{\Phi;\Psi\}$ is a Riesz $\dm$-wavelet in $\dLp{2}$ and $\si(\Phi)=\V_-(\Psi)$.

\item[(iv)] $\V_-(\Psi)$ is shift-invariant and
there exists a subset $\Phi=\{\varphi^1,\ldots,\varphi^r\} \subset \dLp{2}$ such that $\{\Phi;\Psi\}$ is a Riesz $\dm$-wavelet in $\dLp{2}$.

\item[(v)] There exist subsets $\Phi=\{\varphi^1,\ldots,\varphi^r\}$, $\tilde{\Phi}=\{\tilde{\varphi}^1,\ldots,\tilde{\varphi}^r\}$,
    $\tilde{\Psi}=\{\tilde{\psi}^1,\ldots,\tilde{\psi}^s\}$
    of $\dLp{2}$ such that $(\{\tilde{\Phi};\tilde{\Psi}\},\{\Phi;\Psi\})$ is a biorthogonal $\dm$-wavelet in $\dLp{2}$.

\item[(vi)] There exists $\Phi=\{\varphi^1,\ldots,\varphi^r\}\subset \dLp{2}$ such that $\{\Phi;\Psi\}$ is a Riesz $\dm$-wavelet in $\dLp{2}$ and $\si(\Phi)\cup\si(\Psi)\subseteq \si_\dm(\Phi)$.

\item[(vii)] $\Psi$ is a homogeneous Riesz $\dm$-wavelet in $\dLp{2}$, $\V_-(\Psi)$ is shift-invariant, and there exists  $\Phi=\{\phi^1,\ldots,\phi^r\}\subset \dLp{2}$ such that $r\le \frac{s}{|\det(\dm)|-1}$ and $\si(\Phi)\cup \si(\Psi)\subseteq \si_\dm(\Phi)$ (i.e.,
    the refinable structure \eqref{refstr:wavelet} holds with $\phi:=(\phi^1,\ldots,\phi^r)^\tp$ and $\psi:=(\psi^1,\ldots,\psi^s)^\tp$ for some $r\times r$ matrix $\wh{a}$ and $s\times r$ matrix $\wh{b}$ of $2\pi\dZ$-periodic measurable functions on $\dR$).
\end{enumerate}
Moreover, any of the above items (i)--(vii) implies $r_1=r_2=r=\frac{s}{|\det(\dm)|-1}$.
\end{theorem}

\bp (i)$\imply$(ii). The existence of $\tilde{\Psi}$ is guaranteed by Lemma~\ref{lem:V}. Since $(\tilde{\Psi},\Psi)$ is a homogeneous biorthogonal $\dm$-wavelet in $\dLp{2}$,
it was proved in Han~\cite[Theorem~1]{hanmsc94} or \cite[Theorem~1]{han95} (also see \cite{ltw01,lem92} and references therein) that
\be \label{dim:V0}
\dim_{\V_-(\Psi)}(\xi)=\sum_{j=1}^\infty\sum_{\ell=1}^s [\wh{\psi^\ell}((\dm^\tp)^{j}\cdot), \wh{\tilde{\psi}^\ell}((\dm^\tp)^j\cdot)](\xi),\qquad a.e.\,\xi\in \dR.
\ee
%
For completeness, we provide a short self-contained proof here to the identity in \eqref{dim:V0}.
Since $(\tilde{\Psi},\Psi)$ is a homogeneous biorthogonal $\dm$-wavelet in $\dLp{2}$, it must be a homogeneous dual $\dm$-framelet in $\dLp{2}$. As in Proposition~\ref{prop:hdf}, by \cite{rs97,css98,hanbook}, $(\{\tilde{H};\tilde{\Psi}\},\{H;\Psi\})$ is a dual $\dm$-framelet in $\dLp{2}$, where the subsets $H$ and $\tilde{H}$ are defined in \eqref{Psi:H:tH}. Consequently, we have the following representation:
\be \label{general}
f=\sum_{h\in H} \sum_{k\in \dZ}\la f,\tilde{h}(\cdot-k)\ra h(\cdot-k)+\sum_{j=0}^\infty \sum_{\psi\in \Psi} \sum_{k\in \dZ} \la f, \tilde{\psi}_{\dm^j;k}\ra \psi_{\dm^j;k},\qquad f\in \dLp{2}.
\ee
Since $\V_-(\Psi)$ is shift-invariant, we have $\V_-(\Psi)=\si(H)$. By $\V_-(\Psi) =\AS_+(\tilde{\Psi})^\perp$ which is proved in the proof of Lemma~\ref{lem:V}, for every $f\in \V_-(\Psi)$, we must have $\la f, g\ra=0$ for all $g\in \AS_+(\tilde{\Psi})$. Consequently, the representation in \eqref{general} for $f\in \V_-(\Psi)$ becomes
$f=\sum_{h\in H} \sum_{k\in \dZ}\la f,\tilde{h}(\cdot-k)\ra h(\cdot-k)$, which can be rewritten in the frequency domain as
\be \label{special}
\wh{f}(\xi)=\sum_{h\in H} [\wh{f},\wh{\tilde{h}}](\xi) \wh{h}(\xi)=
\sum_{j=1}^\infty \sum_{\ell=1}^s [\wh{f},\wh{\tilde{\psi}^\ell}((\dm^\tp)^j\cdot)](\xi)
\wh{\psi^\ell}((\dm^\tp)^j\xi),\qquad a.e.\,\xi\in \dR
\ee
for $f\in \V_-(\Psi)$, where we used the definition of $H$ and $\tilde{H}$ in \eqref{Psi:H:tH}. Let $\varphi^1,\ldots,\varphi^r\in \dLp{2}$ be given in Proposition~\ref{prop:sis} with $\Phi$ being replaced by $H$. Then $\{\varphi^n(\cdot-k) \setsp k\in \dZ, n=1,\ldots,r\}$ is a (normalized) tight frame for $\si(H)=\V_-(\Psi)$ and $\dim_{\V_-(\Psi)}(\xi)=\sum_{n=1}^r [\wh{\varphi^n},\wh{\varphi^n}](\xi)$.
Since $\varphi^n\in \V_-(\Psi)$, it follows from \eqref{special} that
$\wh{\varphi^n}(\xi)=
\sum_{j=1}^\infty \sum_{\ell=1}^s [\wh{\varphi^n},\wh{\tilde{\psi}^\ell}((\dm^\tp)^j\cdot)](\xi)
\wh{\psi^\ell}((\dm^\tp)^j\xi)$, from which we see that
\[
[\wh{\varphi^n},\wh{\varphi^n}](\xi)
=\sum_{j=1}^\infty \sum_{\ell=1}^s [\wh{\varphi^n},\wh{\tilde{\psi}^\ell}((\dm^\tp)^j\cdot)](\xi)
[\wh{\psi^\ell}((\dm^\tp)^j\cdot), \wh{\varphi^n}](\xi).
\]
Consequently, we deduce from the above identity that
\be \label{dim:eq2}
\dim_{\V_-(\Psi)}(\xi)=\sum_{n=1}^r [\wh{\varphi^n},\wh{\varphi^n}](\xi)=
\sum_{n=1}^r \sum_{j=1}^\infty \sum_{\ell=1}^s [\wh{\varphi^n},\wh{\tilde{\psi}^\ell}((\dm^\tp)^j\cdot)](\xi)
[\wh{\psi^\ell}((\dm^\tp)^j\cdot), \wh{\varphi^n}](\xi).
\ee
On the other hand, since $\{\varphi^n(\cdot-k) \setsp k\in \dZ, n=1,\ldots,r\}$ is a (normalized) tight frame for $\si(H)=\V_-(\Psi)$, by $|\det(\dm)|^{-j} \psi^\ell(\dm^{-j}\cdot) \in \V_-(\Psi)$ for $j\in \N$ and $\ell=1,\ldots,s$,
we have the representation
\[
|\det(\dm)|^{-j} \psi^\ell(\dm^{-j}\cdot)=\sum_{n=1}^r \sum_{k\in \dZ} \la |\det(\dm)|^{-j} \psi^\ell(\dm^{-j}\cdot), \varphi^n(\cdot-k)\ra \varphi^n(\cdot-k).
\]
In the frequency domain, this representation becomes
$\wh{\psi^\ell}((\dm^\tp)^j\xi)=\sum_{n=1}^r [\wh{\psi^\ell}((\dm^\tp)^j\cdot),\wh{\varphi^n}](\xi)
\wh{\varphi^n}(\xi)$. Consequently, we must have
\be \label{dim:eq3}
[\wh{\psi}((\dm^\tp)^j\cdot),\wh{\tilde{\psi}}((\dm^\tp)^j\cdot)](\xi)
=\sum_{n=1}^r [\wh{\psi}((\dm^\tp)^j\cdot),\wh{\varphi^n}](\xi)
[\wh{\varphi^n},\wh{\tilde{\psi}}((\dm^\tp)^j\cdot)](\xi).
\ee
Combining \eqref{dim:eq2} and \eqref{dim:eq3}, we conclude that
\eqref{dim:V0} holds.
This proves (i)$\imply$(ii) with $r_2=r_1$.

(ii)$\imply$(iii). By Lemma~\ref{lem:V}, $\V_-(\Psi)$ is shift-invariant and by \eqref{dim:Psi}, we have $\dim_{\V_-(\Psi)}(\xi)=r_2$ for almost every $\xi\in \dR$.
Consequently, by Proposition~\ref{prop:sis}, there exist $\varphi^1,\ldots,\varphi^r\in \dLp{2}$ with $r=r_2$ such that $S:=\{\varphi^1(\cdot-k),\ldots, \varphi^r(\cdot-k) \setsp k\in \dZ\}$ is an orthonormal basis of $\V_-(\Psi)$. In particular, we have $\si(\Phi)=\V_-(\Psi)$.
Since $\AS_+(\Psi)$ is a Riesz basis of $\V_+(\Psi)$ and $\AS_-(\Psi)$ is a Riesz basis of $\V_-(\Psi)$, we conclude from item (i) of Lemma~\ref{lem:riesz} that $\AS_0(\Phi;\Psi)=S\cup \AS_+(\Psi)$ must be a Riesz basis of $\dLp{2}$.
That is, $\{\Phi;\Psi\}$ must be a Riesz $\dm$-wavelet in $\dLp{2}$ and $\si(\Phi)=\V_-(\Psi)$.

(iii)$\imply$(iv) is trivial.
We now prove (iv)$\imply$(v). Since $\{\Phi;\Psi\}$ is a Riesz $\dm$-wavelet in $\dLp{2}$ and $\dm$ is expansive, by \cite[Theorem~7]{han12}, $\Psi$ must be a homogeneous Riesz $\dm$-wavelet in $\dLp{2}$. Since $\V_-(\Psi)$ is shift-invariant, by Lemma~\ref{lem:V}, there exists $\tilde{\Psi}\subset \dLp{2}$ such that $(\tilde{\Psi},\Psi)$ is a homogeneous biorthogonal $\dm$-wavelet in $\dLp{2}$. Therefore, by Lemma~\ref{lem:V} again, $\V_-(\tilde{\Psi})$ is shift-invariant.
Let $P: \dLp{2}\rightarrow \V_-(\tilde{\Psi})$ be the orthogonal projection and define $\tilde{\Phi}:=P\Phi$.
Since $\V_-(\tilde{\Psi})$ is shift-invariant, $\tilde{S}:=\{\tilde{\varphi}(\cdot-k)
\setsp \tilde{\varphi}\in \tilde{\Phi}\}$ must be the image of
$S:=\{\varphi(\cdot-k)\setsp k\in \dZ,\varphi\in \Phi\}$ under the orthogonal projection operator $P$.
Since $\V_-(\tilde{\Psi})=\V_+(\Psi)^\perp$,
we conclude from item (iii) of Lemma~\ref{lem:riesz} that $\tilde{S}$ must be a Riesz basis of $\V_-(\tilde{\Psi})$. Since $\AS_+(\tilde{\Psi})$ is a Riesz basis of $\V_+(\tilde{\Psi})$ and $\AS_-(\tilde{\Psi})$ is a Riesz basis of $\V_-(\tilde{\Psi})$, by item (i) of Lemma~\ref{lem:riesz}, we conclude that $\AS_0(\tilde{\Phi};\tilde{\Psi})=\tilde{S}\cup \AS_+(\tilde{\Psi})$ is a Riesz basis of $\dLp{2}$.
Since $\AS(\tilde{\Psi})$ and $\AS(\Psi)$ are biorthogonal to each other,
the dual Riesz basis of the Riesz basis $\AS_0(\tilde{\Phi};\tilde{\Psi})$ must be $\AS_0(\Phi;\Psi)$ for some
subset $\Phi\subset \dLp{2}$.
That is, we proved that $(\{\tilde{\Phi};\tilde{\Psi}\},\{\Phi;\Psi\})$ is a biorthogonal $\dm$-wavelet in $\dLp{2}$.


(v)$\imply$(vi) follows directly from \cite[Theorem~8]{han12}.

(vi)$\imply$(vii). Since $\dm$ is expansive, by \cite[Theorem~6]{han12}, $(\tilde{\Psi},\Psi)$ must be a homogeneous biorthogonal $\dm$-wavelet in $\dLp{2}$.
Therefore, by Lemma~\ref{lem:V},
$\V_-(\Psi)$ is shift-invariant and $\Psi$ must be a homogeneous Riesz $\dm$-wavelet in $\dLp{2}$.
By \cite[Theorem~7]{han12} (also, see \cite[Theorem~4.5.1]{hanbook} for dimension one), we must have $s=r(|\det(\dm)|-1)$ and \eqref{refstr:wavelet} holds with $\phi=(\varphi^1,\ldots,\varphi^r)^\tp$.

(vii)$\imply$(i). We only need to prove $\dim_{\V_-(\Psi)}(\xi)=r$ for almost every $\xi\in \dR$. Since $\Psi$ is a homogeneous Riesz $\dm$-wavelet,
we must have $\V_-(\Psi)\cap \si(\Psi)=\{0\}$ and $\len(\si(\Psi))=s$. By the definition of $\V_-(\Psi)$,
it is trivial that $\V_-(\Psi)+\si(\Psi)=\V_\dm(\Psi):=\{
f(\dm\cdot) \setsp f\in \V_-(\Psi)\}$. Therefore, by the definition of the dimension function, we have $\dim_{\V_-(\Psi)}(\xi)+\dim_{\si(\Psi)}(\xi)
=\dim_{\V_\dm(\Psi)}(\xi)$ for almost every $\xi\in \dR$.
By $\dim_{\si(\Psi)}(\xi)=s$ and a similar argument as in \cite{bow03}, we have
\be \label{dim:identity}
\int_{(-\pi,\pi]^d} \dim_{\V_-(\Psi)}(\xi)d\xi+(2\pi)^d s=\int_{(-\pi,\pi]^d} \dim_{\V_\dm(\Psi)}(\xi) d\xi=|\det(\dm)| \int_{(-\pi,\pi]^d} \dim_{\V_-(\Psi)}(\xi)d\xi,
\ee
by $\V_\dm(\Psi)=\{f(\dm\cdot) \setsp f\in \V_-(\Psi)\}$. Since $\V_-(\Psi)\subset \si(\Phi)$ and $\dim_{\si(\Phi)}(\xi)\le r$, we have $\dim_{\V_-(\Psi)}(\xi)\le \dim_{\si(\Phi)}(\xi)\le r$ for almost every $\xi\in \dR$.
We now deduce from \eqref{dim:identity} that
$\int_{(-\pi,\pi]^d} \dim_{\V_-(\Psi)}(\xi)d\xi=\frac{(2\pi)^d s}{|\det(\dm)|-1}$.
By our assumption $r\le \frac{s}{|\det(\dm)|-1}$ and  $\dim_{\V_-(\Psi)}(\xi)\le r$, we have
\[
(2\pi)^d r \le \frac{(2\pi)^d s}{|\det(\dm)|-1}
=
\int_{(-\pi,\pi]^d} \dim_{\V_-(\Psi)}(\xi)d\xi\le
\int_{(-\pi,\pi]^d} r d\xi= (2\pi)^d r.
\]
Since $\dim_{\V_-(\Psi)}(\xi)\le r$ for almost every $\xi\in \R$, we conclude from the above inequalities that $\dim_{\V_-(\Psi)}(\xi)=r$ for almost every $\xi\in \dR$.
\ep

For a homogeneous biorthogonal $\dm$-wavelet $(\tilde{\Psi},\Psi)$ in $\dLp{2}$, Theorem~\ref{thm:hw} tells us that there exist $\Phi,\tilde{\Phi}\subset \dLp{2}$ such that $(\{\tilde{\Phi};\tilde{\Psi}\},\{\Phi;\Psi\})$ is a biorthogonal $\dm$-wavelet in $\dLp{2}$ if and only if the identity \eqref{dim:Psi} holds. For a homogeneous Riesz $\dm$-wavelet $\Psi$ in $\dLp{2}$, under the extra condition that $\V_-(\Psi)$ is shift-invariant, Theorem~\ref{thm:hw} tells us that there exists $\Phi\subset \dLp{2}$ such that $\{\Phi;\Psi\}$ is a Riesz $\dm$-wavelet in $\dLp{2}$ if and only if \eqref{dim:Psi} holds. Without the assumption that $\V_-(\Psi)$ is shift-invariant, for a given homogeneous Riesz $\dm$-wavelet $\Psi$ in $\dLp{2}$,
it remains unclear to us what is a necessary and sufficient condition for the existence of a Riesz $\dm$-wavelet $\{\Phi;\Psi\}$ in $\dLp{2}$.

One of the main motivations for linking homogeneous wavelets to nonhomogeneous wavelets is the following known result (see \cite[Theorem~4.5.1]{hanbook} and \cite[Theorem~7]{han12}), showing that a nonhomogeneous wavelet intrinsically has the refinable structure.

\begin{theorem}\label{thm:nbw}
Let $\dm$ be a $d\times d$ invertible integer matrix. Let $\Phi=\{\phi^1,\ldots,\phi^r\}$, $\tilde{\Phi}:=\{\tilde{\phi}^1,\ldots,\tilde{\phi}^r\}$,
$\Psi=\{\psi^1,\ldots,\psi^s\}$ and $\tilde{\Psi}=\{\tilde{\psi}^1,\ldots,\tilde{\psi}^s\}$ be finite subsets of $\dLp{2}$.
Define    $\phi:=(\phi^1,\ldots,\phi^r)^\tp$, $\tilde{\phi}:=(\tilde{\phi}^1,\ldots,\tilde{\phi}^r)^\tp$ and
$\psi:=(\psi^1,\ldots,\psi^s)^\tp$, $\tilde{\psi}:=(\tilde{\psi}^1,\ldots,\tilde{\psi}^s)^\tp$. Then $(\{\tilde{\Phi};\tilde{\Psi}\},\{\Phi;\Psi\})$
is a biorthogonal $\dm$-wavelet in $\dLp{2}$ if and only if
\begin{enumerate}
\item[(1)] $\lim_{j\to \infty}
\la \wh{\tilde{\phi}}((\dm^\tp)^{-j}\cdot)^\tp
\ol{\wh{\phi}((\dm^\tp)^{-j}\cdot)}, h\ra=\la 1,h\ra$ for all compactly supported $\mathscr{C}^\infty$ functions $h$ on $\dR$;
\item[(2)] The vector functions $\phi$ and $\tilde{\phi}$ are biorthogonal to each other: $\la \tilde{\phi},\phi(\cdot-k)\ra:=\int_{\dR}
    \tilde{\phi}(x)\ol{\phi(x-k)}^\tp dx$ $=\td(k)I_r$ for all $k\in \dZ$, where $\td(0)=1$ and $\td(k)=0$ for all $k\ne 0$;
\item[(3)] There exist $r\times r$ matrices $\wh{a}, \wh{\tilde{a}}$ and $s\times r$ matrices $\wh{b},\wh{\tilde{b}}$ of $2\pi$-periodic measurable functions in $\dTLp{2}$ such that \eqref{reffunc} and \eqref{filterbank} are satisfied,
%
%
and $(\{\wh{\tilde{a}};\wh{\tilde{b}}\},\{\wh{a};\wh{b}\})$
is a biorthogonal $\dm$-wavelet filter bank, i.e., $s=r(|\det(\dm)|-1)$ and
$\wh{\tilde{a}}(\xi)^\tp\ol{\wh{a}(\xi+2\pi \omega)}+\wh{\tilde{b}}(\xi)^\tp \ol{\wh{b}(\xi+2\pi \omega)}=\td(\omega)I_r$
a.e. $\xi\in \dR$
for all $\omega\in \dmfc:=[(\dm^\tp)^{-1}\dZ]\cap[0,1)^d$;

\item[(4)] $\AS_0(\Phi;\Psi)$ and $\AS_0(\tilde{\Phi};\tilde{\Psi})$ are Bessel sequences in $\dLp{2}$.
\end{enumerate}
\end{theorem}

Theorem~\ref{thm:nbw} shows that
a nonhomogeneous biorthogonal wavelet has
the intrinsic refinable structure. As a consequence, under the condition that $\V_-(\Psi)$ is shift-invariant, it is not surprising that the statements in items (i), (ii) and (vii) of Theorem~\ref{thm:hw} are very similar to the results in \cite{bow03,han95,hanmsc94,hwbook96,kkl01,lem92,zal99}
linking homogeneous Riesz wavelets and homogeneous biorthogonal wavelets with a (generalized) multiresolution analysis.
However, a nonhomogeneous wavelet in $\dLp{2}$ is not always associated with a generalized multiresolution analysis and therefore, it remains unclear whether the condition that $\V_-(\Psi)$ is shift-invariant in items (i), (iv) and (vii) of Theorem~\ref{thm:hw} can be dropped or not to obtain a stronger conclusion for linking homogeneous wavelets to nonhomogeneous ones.

For the special case of homogeneous orthogonal $\dm$-wavelets in Theorem~\ref{thm:hw}, the space $\V_-(\Psi)$ is always automatically shift-invariant and hence,
the condition on $\V_-(\Psi)$ can be dropped. More precisely, we have the following result.

\begin{cor}\label{cor:ow}
Let $\dm$ be a $d\times d$ expansive integer matrix and $\Psi=\{\psi^1,\ldots,\psi^s\}$ be a finite subset of $\dLp{2}$. Then the following statements are equivalent:
\begin{enumerate}
\item[(1)] $\Psi$ is a homogeneous orthogonal $\dm$-wavelet in $\dLp{2}$ and for some constant $r_1$,
\be \label{dim:Psi:ow}
\sum_{j=1}^\infty\sum_{\ell=1}^s [\wh{\psi^\ell}((\dm^\tp)^{j}\cdot), \wh{\psi^\ell}((\dm^\tp)^j\cdot)](\xi)=r_1,\quad a.e.\, \xi\in \dR.
\ee

\item[(2)] There exists a subset $\Phi=\{\varphi^1,\ldots,\varphi^r\} \subset \dLp{2}$ such that $\{\Phi;\Psi\}$ is an orthogonal $\dm$-wavelet in $\dLp{2}$.

\item[(3)] $\AS_+(\Psi):=\{\psi_{\dm^j;k} \setsp j\in \N\cup\{0\},k\in \dZ,\psi\in \Psi\}$ is an orthonormal system in $\dLp{2}$ and there exists $\Phi=\{\varphi^1,\ldots,\varphi^r\}\subset \dLp{2}$ such that $\{\Phi;\Psi\}$ is a Riesz $\dm$-wavelet in $\dLp{2}$.

\item[(4)] $\Psi$ is a homogeneous orthogonal $\dm$-wavelet in $\dLp{2}$ and there exists $\Phi=\{\phi^1,\ldots,\phi^r\}\subset \dLp{2}$ such that
    $r\le \frac{s}{|\det(\dm)|-1}$ and $\si(\Phi)\cup\si(\Psi)\subseteq \si_\dm(\Phi)$ (i.e.,
    the refinable structure in \eqref{refstr:wavelet} holds with $\phi:=(\phi^1,\ldots,\phi^r)^\tp$ and $\psi:=(\psi^1,\ldots,\psi^s)^\tp$ for some $r\times r$ matrix $\wh{a}$ and $s\times r$ matrix $\wh{b}$ of $2\pi\dZ$-periodic measurable functions on $\dR$).
\end{enumerate}
Moreover, if any of the above items (1)--(4) is satisfied, then we must have $r_1=r=\frac{s}{|\det(\dm)|-1}$.
\end{cor}

\bp Since $\Psi$ is a homogeneous orthogonal $\dm$-wavelet in $\dLp{2}$, by Lemma~\ref{lem:V} or the identity $\V_-(\Psi)=\V_+(\Psi)^\perp$,
we see that $\V_-(\Psi)$ must be shift-invariant. Moreover, if $\{\Phi;\Psi\}$ is an orthogonal $\dm$-wavelet in $\dLp{2}$, then it is also trivial that $\si(\Phi)=\V_-(\Psi)$.
Now all the claims follow directly from Theorem~\ref{thm:nbw}.
\ep

Let $\Psi\subset \dLp{2}$ be a homogeneous orthogonal $\dm$-wavelet in $\dLp{2}$ and define a subset $H:=\{ |\det(\dm)|^{-j} \psi(\dm^{-j}\cdot) \setsp j\in \N, \psi\in \Psi\}$ as in \eqref{Psi:H}. Then $\Psi$ must be a homogeneous tight $\dm$-framelet in $\dLp{2}$. Suppose that $r=\len(\si(H))<\infty$. By Theorem~\ref{thm:htf}, there exists a subset $\Phi\subset \si(H)$ with $\#\Phi\le r$ such that $\{\Phi;\Psi\}$ is a tight $\dm$-framelet in $\dLp{2}$. On the other hand, according to the characterization of nonhomogeneous tight framelets in Theorem~\ref{thm:ntf}, the sets $\Phi$ and $\Psi$ must have the refinable structure $\si(\Phi)\cup \si(\Psi)\subseteq \si_\dm(\Phi)$. As a consequence, every homogeneous orthogonal $\dm$-wavelet $\Psi$ in $\dLp{2}$ satisfying $\len(\si(H))<\infty$ (this condition can be dropped if we allow $\Phi$ to be an infinitely countable set) must be derived from the refinable structure through a generalized tight $\dm$-framelet filter bank in item (iii) of Theorem~\ref{thm:ntf}.
However, $\{\Phi;\Psi\}$ may not be an orthogonal $\dm$-wavelet in $\dLp{2}$ and $s>r(\#\Psi)-r$ may happen. According to Corollary~\ref{cor:ow},
there exists a subset $\Phi\subset \dLp{2}$ such that $\{\Phi;\Psi\}$ is a (nonhomogeneous) orthogonal $\dm$-wavelet in $\dLp{2}$ if and only if \eqref{dim:Psi:ow} is satisfied with the constant integer $r_1=\frac{s}{|\det(\dm)|-1}$.

For $\dm=2$ and $\psi\in \Lp{2}$, if $\{\psi\}$ is a homogeneous orthogonal $2$-wavelet in $\Lp{2}$ such that $\wh{\psi}$ is continuous and satisfies
$|\wh{\psi}(\xi)|=\bo(|\xi|^\gep)$ and
$\wh{\psi}(\xi)|=\bo((1+|\xi|)^{-\gep-1/2})$
for some $\gep>0$,
then it is known in \cite{lem92} that $\sum_{j=1}^\infty [\wh{\psi}(2^j\cdot),\wh{\psi}(2^j\cdot)](\xi)=1$, i.e., \eqref{dim:Psi:ow} holds for all $\xi\in \R\bs (2\pi \Z)$ (see \cite{aus95} for its generalization to homogeneous biorthogonal $2$-wavelets).
Hence, under such extra conditions, there always exists a function $\phi\in \Lp{2}$ such that $\{\{\phi\}; \{\psi\}\}$ is an orthogonal $2$-wavelet in $\Lp{2}$.

We finish this paper by presenting an example to demonstrate that not every homogeneous Riesz wavelet $\Psi$ in $\dLp{2}$ has an associated nonhomogeneous wavelet $\{\Phi;\Psi\}$ in $\dLp{2}$.

\noindent \textbf{Example.}
Define $\wh{\psi}:=\chi_{K\cup(-K)}$ with $K:=[\frac{4\pi}{7},\pi]\cup[4\pi,\frac{32\pi}{7}]$.
Then $\psi$ is the known Journ\'e homogeneous orthogonal $2$-wavelet and $\psi$ cannot be generated by a multiresolution analysis (see \cite{daubook92,meybook90}).
Note that $\V_-(\psi)$ is shift-invariant.
By calculation, it is known in \cite[Example~3.6]{bjmp05} that
\[
\dim_{\V_-(\psi)}(\xi)=
\chi_{[-\pi,-\frac{6\pi}{7}]\cup [-\frac{4\pi}{7},\frac{4\pi}{7}]
\cup[\frac{6\pi}{7},\pi]}(\xi)+
2\chi_{[-\frac{2\pi}{7},\frac{2\pi}{7}]}(\xi), \qquad \xi\in [-\pi,\pi)
\]
except possibly at the endpoints of the intervals. Therefore, $\len(\V_-(\psi))=\|\dim_{\V_-(\psi)}(\cdot)\|_{\Lp{\infty}}=2$ and by Theorem~\ref{thm:htf}, there exist $\phi^1,\phi^2\in \V_-(\psi)$ such that $\{\Phi; \{\psi\}\}$ with $\Phi:=\{\phi^1,\phi^2\}$ is a tight $2$-framelet in $\Lp{2}$. Moreover, by Theorem~\ref{thm:ntf}, the tight $2$-framelet $\{\Phi,\{\psi\}\}$ must have the refinable structure $\si(\Phi)\cup \si(\psi)\subseteq \si_2(\Phi)$. That is,
$\wh{\phi}(2\xi)=\wh{a}(\xi)\wh{\phi}(\xi)$
and $\wh{\psi}(2\xi)=\wh{b}(\xi)\wh{\phi}(\xi)$ are satisfied for some $2\times 2$ matrix $\wh{a}$ and some $1 \times 2$ vector $\wh{b}$ of $2\pi$-periodic measurable functions on $\R$, where $\phi:=(\phi^1,\phi^2)^\tp$.
Moreover, $\{\wh{a};\wh{b}\}$ is a generalized tight $2$-framelet filter bank (see item (iii) of Theorem~\ref{thm:ntf} for definition).
But by Corollary~\ref{cor:ow},
since $\dim_{\V_-(\psi)}$ is not a constant function (and hence \eqref{dim:Psi:ow} fails), there does not exist a subset $\Phi\subset \Lp{2}$ such that $\{\Phi;\{\psi\}\}$ is a Riesz $2$-wavelet in $\Lp{2}$.

\end{document}